\documentclass [12pt] {amsart}
\usepackage{latexsym,amsfonts,amssymb,euler,euscript,amscd,epsfig,isolatin1,
supertabular,rotating,epic}
\usepackage [all]{xy}
\usepackage{fullpage}
\theoremstyle{plain}
\newtheorem {theorem}{Theorem}[subsection]
\newtheorem {Thm}[theorem] {Theorem}
\newtheorem {Lem}[theorem] {Lemma}
\newtheorem {Prop}[theorem] {Proposition}
\newtheorem {Cor}[theorem] {Corollary}
\theoremstyle {definition}
\newtheorem {Def}[theorem] {Definition}
\newtheorem {Not}[theorem] {Notation}
\newtheorem {Con}[theorem] {Construction}
\newtheorem {Rem}[theorem] {Remark}
\newtheorem {Exa}[theorem] {Example}
\newenvironment{Pf}[1]{{\sc Proof#1:}}{\qed\medskip\\}

\newcommand {\Hom} {\operatorname {Hom}}
\newcommand {\equalizer} {\operatorname{eq}}
\newcommand {\coeq} {\operatorname{coeq}}
\newcommand {\im}{\operatorname{im}}
\newcommand {\hgt} {\operatorname{lgt}}

\def\HoRKan {\mathop{HoRKan}}
\newcommand {\horkan}[1]{\HoRKan_{#1}}
\def\HoLKan {\mathop{\operatorname{HoLKan}}}
\newcommand {\holkan}[1]{\HoLKan_{#1}}

\def\LKan {\mathop{\operatorname{LKan}}}
\newcommand {\lkan} [1] {\LKan_{#1}}
\def\HoPfeil{{\unitlength 1em\begin{picture}(0,.1)
\put(0,.1){\vector(1,0){2.65}}\end{picture}}}
\def\Hocolimname{{\unitlength.1em
\raisebox{-2.7\unitlength}{\begin{picture}(27,9.5)(0,0)
\put(0,2.7){$\operatorname{Holim}$}
\put(.1,-.1){\HoPfeil} 
\end{picture}}}}
\def\Hocolim{\mathop{\Hocolimname}}
\def\Pfeil{{\unitlength 1em\begin{picture}(0,.1)
\put(0,.1){\vector(1,0){1.5}}\end{picture}}}
\def\colimname{{\unitlength.1em
\raisebox{-2.7\unitlength}{\begin{picture}(15.5,9.5)(0,0)
\put(0,2.7){$\operatorname{lim}$}
\put(.05,-.1){\Pfeil} 
\end{picture}}}}
\def\Colim{\mathop{\colimname}}
\newcommand {\colim}[1] {\Colim_{#1}}
\newcommand {\colimlim}[1] {\Colim\limits_{#1}}
\newcommand {\hocolim}[1] {\Hocolim_{#1}}

\newcommand {\cone} {\operatorname{cone}}
\newcommand {\Cone} {\operatorname{\mathfrak{Cone}}}
\def\SqCp{{\unitlength.1em\begin{picture}(10,8)
\put(0,0){\rule{1em}{.06em}}                     
\put(2,0){\rule{.06em}{.8em}}                   
\put(7.4,0){\rule{.06em}{.8em}}                  
\end{picture}}}

\def\sQcup{\mathrel{\SqCp}}
\newcommand{\push}[1] {\underset {#1} \sQcup}
\newcommand {\diag} {\operatorname{diag}}

\newcommand {\E} {\operatorname E(1)_*}
\newcommand {\EE} {\operatorname E(1)_*E(1)}
\newcommand {\elok} {category of $E(1)$-local spectra } 
\newcommand {\CA}[1] {(C_N\times C_N)\to #1}

\newcommand {\CN} {{\mathcal C}^{2p-2}(\eeO)}
  
\newcommand {\SE} {\mathcal S_{E(1)}}
\newcommand {\DD} {\mathcal{D}} 
\newcommand {\MM} {\mathcal M}
\newcommand {\EM} {\mathcal M_{S_{E(1)}}}

\newcommand {\eeO} {\operatorname{Comod}^{0}_{\EE}}
\newcommand {\eee} {\operatorname{Comod}_{\EE}}

\newcommand {\C} {C^\bullet}
\newcommand {\Z} {Z^\bullet}
\newcommand {\B} {B^\bullet}

\renewcommand {\b} {{\mathcal B}^{s+1}}
\newcommand {\z} {{\mathcal Z}^s}
\newcommand {\sBB} {{\mathcal B}^{s+1}\wedge\tilde{\mathcal B}^{t+1}}
\newcommand {\szz} {{\mathcal Z}^s\wedge\tilde{\mathcal Z}^t}

\renewcommand {\c} {{\mathcal C}^s}
\newcommand {\scc} {{\mathcal C}^s\wedge\tilde{\mathcal C}^t}

\newcommand {\sbc} {{\mathcal B}^s\wedge\tilde{\mathcal C}^t}
\newcommand {\scb} {{\mathcal C}^s\wedge\tilde{\mathcal B}^t}
\newcommand {\bzbbzb} {\b\wedge\tilde{\mathcal Z}^t\push{\b\wedge
                       \tilde{\mathcal B}^{t+1}}
                       \z\wedge\tilde{\mathcal B}^{t+1}}
\newcommand {\sBz} {{\mathcal B}^{s+1}\wedge\tilde{\mathcal Z}^t}
\newcommand {\szB} {{\mathcal Z}^s\wedge\tilde{\mathcal B}^{t+1}}
\newcommand {\sbz} {{\mathcal B}^s\wedge\tilde{\mathcal Z}^t}
\newcommand {\szb} {{\mathcal Z}^s\wedge\tilde{\mathcal B}^t}
\newcommand {\sbb} {{\mathcal B}^s\wedge\tilde{\mathcal B}^t}
\newcommand {\sBb} {{\mathcal B}^{s+1}\wedge\tilde{\mathcal B}^t}
\newcommand {\sbB} {{\mathcal B}^s\wedge\tilde{\mathcal B}^{t+1}}
\newcommand {\PP} {\tilde P}
\newcommand {\AAA} {(A\wedge\tilde A)}
\newcommand {\CCto} {C_N\times C_N \to} 
\newcommand {\Ob} {\operatorname{Ob}}
\newcommand {\Mor} {\operatorname{Mor}} 

\newcommand {\Rek} {\mathcal{R}ec}
\newcommand{\Ext} {\operatorname{Ext}}
\newcommand {\Ho} {\operatorname{Ho}}

\newcommand {\specialmap} [4] {\text {$ #1\negmedspace : #2 #3 #4 $}}
\newcommand {\map} [3] {\specialmap {#1} {#2}{\longrightarrow} {#3}}

\newcommand {\id}{\operatorname{id}}

\newcommand {\zn} {\bf Z_{\!/\!\raisebox{-.2ex}{$_{2p-2}$}}}

\newcommand {\at} [1] {\arrowvert_{#1}}

\newcommand{\Vv}{{
\unitlength .8ex
\begin{picture}(3,0.2)
\def\Kr{\circle{.4}}
\put(0.5,1){\Kr}
\put(2.5,1){\Kr}
\put(1.5,0){\Kr}
\drawline(1.3586,0.1414)(0.6414,0.8586)
\drawline(1.6414,0.1414)(2.3586,0.8586)
\end{picture}}}

\newcommand{\halfVv}{{
\unitlength .8ex
\begin{picture}(2,0.2)
\def\Kr{\circle{.4}}
\put(1.5,1){\Kr}
\put(0.5,0){\Kr}
\drawline(0.6414,0.1414)(1.3586,0.8586)
\end{picture}}}

\newcommand{\VVv}{{
\unitlength .7em
\begin{picture}(2,0.2)
\def\Kr{\circle{.4}}
\put(0,1){\Kr}
\put(2,1){\Kr}
\put(1,0){\Kr}
\drawline(0.8586,0.1414)(0.1414,0.8586)
\drawline(1.1414,0.1414)(1.8586,0.8586)
\end{picture}}}

\newcommand{\ww}{{
\unitlength .8ex
\begin{picture}(4,0.2)
\def\Kr{\circle{.4}}
\put(0,1){\Kr}
\put(2,1){\Kr}
\put(4,1){\Kr}
\put(1,0){\Kr}
\put(3,0){\Kr}
\drawline(0.8586,0.1414)(0.1414,0.8586)
\drawline(1.1414,0.1414)(1.8586,0.8586)
\drawline(2.8586,0.1414)(2.1414,0.8586)
\drawline(3.1414,0.1414)(3.8586,0.8586)
\end{picture}}}

\newcommand{\Ww}{{
\unitlength .7em
\begin{picture}(4,0.2)
\def\Kr{\circle{.4}}
\put(0,1){\Kr}
\put(2,1){\Kr}
\put(4,1){\Kr}
\put(1,0){\Kr}
\put(3,0){\Kr}
\drawline(0.8586,0.1414)(0.1414,0.8586)
\drawline(1.1414,0.1414)(1.8586,0.8586)
\drawline(2.8586,0.1414)(2.1414,0.8586)
\drawline(3.1414,0.1414)(3.8586,0.8586)
\end{picture}}}

\newcommand {\jw} {j_\ww}

\newcommand{\Mv}{\MM^\Vv}

\newcommand{\Vst}{{
\unitlength .6em
\begin{picture}(3,0.8)
\def\Kr{\circle{.4}}
\put(0.2,1){$\star$}
\put(2.5,1){\Kr}
\put(1.5,0){\Kr}
\drawline(1.3586,0.1414)(0.6414,0.8586)
\drawline(1.6414,0.1414)(2.3586,0.8586)
\end{picture}}}

\newcommand{\VC}{
{\unitlength .8ex
\begin{picture}(2,2)
\def\Kr{\circle{.4}}
\put(0,0.5){\Kr}
\put(2,0.5){\Kr}
\put(1,0){\Kr}
\put(1,1.5){\Kr}
\drawline(0.8211,0.0894)(0.1789,0.4106)
\drawline(1.1789,0.0894)(1.8211,0.4106)
\drawline(0.1414,0.6414)(0.8586,1.3586)
\drawline(1.8586,0.6414)(1.1414,1.3586)
\end{picture}}}

\newcommand{\Vo}{
{\unitlength .8ex
\begin{picture}(5,2)
\def\Kr{\circle{.4}}
\put(1.5,0.5){\Kr}
\put(3.5,0.5){\Kr}
\put(2.5,0){\Kr}
\put(2.5,1.5){\Kr}
\put(0.5,0){\Kr}
\put(4.5,0){\Kr}
\drawline(2.3211,0.0894)(1.6789,0.4106)
\drawline(2.6789,0.0894)(3.3211,0.4106)
\drawline(1.6414,0.6414)(2.3586,1.3586)
\drawline(3.3586,0.6414)(2.6414,1.3586)
\drawline(0.6789,0.0894)(1.3211,0.4106)
\drawline(4.3211,0.0894)(3.6789,0.4106)
\drawline(4.6789,0.0894)(4.9,0.2)
\drawline(0.3211,0.0894)(0.1,0.2)
\end{picture}}}

\newcommand{\lo}{l_\Vo}
\newcommand{\jo}{j_\Vo}
\newcommand{\po}{p_\Vo}

\newcommand{\VO}{
{\unitlength .7em
\begin{picture}(5,2)
\def\Kr{\circle{.4}}
\put(1.5,0.5){\Kr}
\put(3.5,0.5){\Kr}
\put(2.5,0){\Kr}
\put(2.5,1.5){\Kr}
\put(0.5,0){\Kr}
\put(4.5,0){\Kr}
\drawline(2.3211,0.0894)(1.6789,0.4106)
\drawline(2.6789,0.0894)(3.3211,0.4106)
\drawline(1.6414,0.6414)(2.3586,1.3586)
\drawline(3.3586,0.6414)(2.6414,1.3586)
\drawline(0.6789,0.0894)(1.3211,0.4106)
\drawline(4.3211,0.0894)(3.6789,0.4106)
\drawline(4.6789,0.0894)(4.9,0.2)
\drawline(0.3211,0.0894)(0.1,0.2)
\end{picture}}}

\newcommand{\Vy}{
{\unitlength .8ex
\begin{picture}(3,2)
\def\Kr{\circle{.4}}
\put(1.5,0.5){\Kr}
\put(1.5,1.5){\Kr}
\put(2.5,0){\Kr}
\put(0.5,0){\Kr}
\drawline(1.5,0.7)(1.5,1.3)
\drawline(0.6789,0.0894)(1.3211,0.4106)
\drawline(2.3211,0.0894)(1.6789,0.4106)
\drawline(2.6789,0.0894)(2.9,0.2)
\drawline(0.3211,0.0894)(0.1,0.2)
\end{picture}}}

\newcommand{\VY}{
{\unitlength .7em
\begin{picture}(3,2)
\def\Kr{\circle{.4}}
\put(1.5,0.5){\Kr}
\put(1.5,1.5){\Kr}
\put(2.5,0){\Kr}
\put(0.5,0){\Kr}
\drawline(1.5,0.7)(1.5,1.3)
\drawline(0.6789,0.0894)(1.3211,0.4106)
\drawline(2.3211,0.0894)(1.6789,0.4106)
\drawline(2.6789,0.0894)(2.9,0.2)
\drawline(0.3211,0.0894)(0.1,0.2)
\end{picture}}}

\newcommand{\py}{p_\Vy}

\newcommand{\pC}{p_\VC}
\newcommand{\iC}{i_\VC}
\newcommand{\lC}{l_\VC}
\newcommand{\jC}{j_\VC}

\newcommand{\II}{{
\unitlength .7em
\begin{picture}(1,1)
\def\Kr{\circle{.4}}
\put(0.5,1){\Kr}
\put(0.5,0){\Kr}
\drawline(0.5,0.2)(0.5,0.8)
\end{picture}}}

\newcommand{\kegink}{{
\unitlength .8em
\begin{picture}(2,1)
\def\Kr{\circle{.4}}
\put(0.25,.75){$\star$}
\put(0.5,0){\Kr}
\drawline(0.5,0.2)(0.5,0.8)
\put(1.5,1){\Kr}
\put(1.5,0){\Kr}
\drawline(1.5,0.2)(1.5,0.8)
\drawline(.7,0)(1.3,0)
\drawline(.7,1)(1.3,1)
\end{picture}}}

\newcommand{\kegabb}{{
\unitlength .8em
\begin{picture}(2,1)
\def\Kr{\circle{.4}}
\put(0.5,1){\Kr}
\put(0.5,0){\Kr}
\drawline(0.5,0.2)(0.5,0.8)
\put(1.5,1){\Kr}
\put(1.25,-.25){$\star$}
\drawline(1.5,0.2)(1.5,0.8)
\drawline(.7,0)(1.3,0)
\drawline(.7,1)(1.3,1)
\end{picture}}}

\newcommand{\Ii}{{
\unitlength .8ex
\begin{picture}(1,1)
\def\Kr{\circle{.4}}
\put(0.5,1.2){\Kr}
\put(.5,0.2){\Kr}
\drawline(0.5,0.4)(0.5,1)
\end{picture}}}
  
\newcommand{\MI}{\MM^\Ii}

\newcommand{\III}[2] {\Setunitlength\raisebox{2\unitlength}
{\xymatrix{
{#1}\ar@{-}[d]\\{#2}}
}}

\newcommand{\VV}[3]
{\Setunitlength\raisebox{2\unitlength}
{\xymatrix{
#1 && #2 \\
&#3\ar@{-}[lu]\ar@{-}[ru]
}}} 

\newcommand{\butterfly}[9]
{
{\Setunitlength
\begin{picture}(26,13)
\def\Kr{\circle{.4}}
\put(13,3){\Kr} 
\put(23,4){\Kr}
\put(13,5){\Kr}
\put(13,11){\Kr}
\put(11,8){\Kr}
\put(15,8){\Kr}
\put(5,7){\Kr}
\put(3,4){\Kr}
\put(21,7){\Kr}
\put(11.1468,8.2202){\line(2,3){1.7065}}
\put(13.1468,5.2202){\line(2,3){1.7065}}
\put(14.8532,8.2202){\line(-2,3){1.7065}}
\put(12.8532,5.2202){\line(-2,3){1.7065}}
\put(3.1468,4.2202){\line(2,3){1.7065}}
\put(22.8532,4.2202){\line(-2,3){1.7065}}
\put(13.1788,3.0894){\line(2,1){7.6422}}
\put(3.1788,4.0894){\line(2,1){7.6422}}
\put(5.1788,7.0894){\line(2,1){7.6422}}
\put(12.8211,3.0894){\line(-2,1){7.6422}}
\put(22.8211,4.0894){\line(-2,1){7.6422}}
\put(20.8211,7.0894){\line(-2,1){7.6422}}
\put(23.1468,4.2202){\line(2,3){1}}
\put(23.1788,4.0894){\line(2,1){1}}
\put(2.8532,4.2202){\line(-2,3){1}}
\put(2.8211,4.0894){\line(-2,1){1}}
\put(1.5,6){\makebox(0,0)[rb]{$\cdots$}}
\put(24.5,6){\makebox(0,0)[lb]{$\cdots$}}
\put(13,11.8){\makebox(0,0)[cb]{$#1$}}
\put(13,2){\makebox(0,0)[ca]{$#9$}}
\put(15.8,8){\makebox(0,0)[lb]{$#3$}}
\put(10.2,8){\makebox(0,0)[rb]{$#2$}}
\put(21.8,7){\makebox(0,0)[lb]{$#6$}}
\put(4.2,7){\makebox(0,0)[rb]{$#5$}}
\put(23.8,3){\makebox(0,0)[la]{$#8$}}
\put(2.2,3){\makebox(0,0)[ra]{$#7$}}
\put(14,5){\makebox(0,0)[la]{$#4$}} 
\end{picture}}
}   

\newcommand{\npluseinsbutterfly}[9]
{
{\Setunitlength
\begin{picture}(26,13)
\def\Kr{\circle{.4}}
\def\ID{\circle*{.4}}
\put(13,3){\ID} 
\put(23,4){\ID}
\put(13,5){\Kr}
\put(13,11){\Kr}
\put(11,8){\Kr}
\put(15,8){\Kr}
\put(5,7){\ID}
\put(3,4){\ID}
\put(21,7){\ID}
\put(11.1468,8.2202){\line(2,3){1.7065}}
\put(13.1468,5.2202){\line(2,3){1.7065}}
\put(14.8532,8.2202){\line(-2,3){1.7065}}
\put(12.8532,5.2202){\line(-2,3){1.7065}}
\put(3.1468,4.2202){\line(2,3){1.7065}}
\put(22.8532,4.2202){\line(-2,3){1.7065}}
\put(13.1788,3.0894){\line(2,1){7.6422}}
\put(3.1788,4.0894){\line(2,1){7.6422}}
\put(5.1788,7.0894){\line(2,1){7.6422}}
\put(12.8211,3.0894){\line(-2,1){7.6422}}
\put(22.8211,4.0894){\line(-2,1){7.6422}}
\put(20.8211,7.0894){\line(-2,1){7.6422}}
\put(23.1468,4.2202){\line(2,3){1}}
\put(23.1788,4.0894){\line(2,1){1}}
\put(2.8532,4.2202){\line(-2,3){1}}
\put(2.8211,4.0894){\line(-2,1){1}}
\put(1.5,6){\makebox(0,0)[rb]{$\cdots$}}
\put(24.5,6){\makebox(0,0)[lb]{$\cdots$}}
\put(13,11.8){\makebox(0,0)[cb]{$#1$}}
\put(13,2){\makebox(0,0)[ca]{$#9$}}
\put(15.8,8){\makebox(0,0)[lb]{$#3$}}
\put(10.2,8){\makebox(0,0)[rb]{$#2$}}
\put(21.8,7){\makebox(0,0)[lb]{$#6$}}
\put(4.2,7){\makebox(0,0)[rb]{$#5$}}
\put(23.8,3){\makebox(0,0)[la]{$#8$}}
\put(2.2,3){\makebox(0,0)[ra]{$#7$}}
\put(14,5){\makebox(0,0)[la]{$#4$}} 
\end{picture}}
}   

\newcommand{\ganzschwarzbutterfly}[9]
{
{\Setunitlength
\begin{picture}(26,13)
\def\ID{\circle*{.4}}
\put(13,3){\ID} 
\put(23,4){\ID}
\put(13,5){\ID}
\put(13,11){\ID}
\put(11,8){\ID}
\put(15,8){\ID}
\put(5,7){\ID}
\put(3,4){\ID}
\put(21,7){\ID}
\put(11.1468,8.2202){\line(2,3){1.7065}}
\put(13.1468,5.2202){\line(2,3){1.7065}}
\put(14.8532,8.2202){\line(-2,3){1.7065}}
\put(12.8532,5.2202){\line(-2,3){1.7065}}
\put(3.1468,4.2202){\line(2,3){1.7065}}
\put(22.8532,4.2202){\line(-2,3){1.7065}}
\put(13.1788,3.0894){\line(2,1){7.6422}}
\put(3.1788,4.0894){\line(2,1){7.6422}}
\put(5.1788,7.0894){\line(2,1){7.6422}}
\put(12.8211,3.0894){\line(-2,1){7.6422}}
\put(22.8211,4.0894){\line(-2,1){7.6422}}
\put(20.8211,7.0894){\line(-2,1){7.6422}}
\put(23.1468,4.2202){\line(2,3){1}}
\put(23.1788,4.0894){\line(2,1){1}}
\put(2.8532,4.2202){\line(-2,3){1}}
\put(2.8211,4.0894){\line(-2,1){1}}
\put(1.5,6){\makebox(0,0)[rb]{$\cdots$}}
\put(24.5,6){\makebox(0,0)[lb]{$\cdots$}}
\put(13,11.8){\makebox(0,0)[cb]{$#1$}}
\put(13,2){\makebox(0,0)[ca]{$#9$}}
\put(15.8,8){\makebox(0,0)[lb]{$#3$}}
\put(10.2,8){\makebox(0,0)[rb]{$#2$}}
\put(21.8,7){\makebox(0,0)[lb]{$#6$}}
\put(4.2,7){\makebox(0,0)[rb]{$#5$}}
\put(23.8,3){\makebox(0,0)[la]{$#8$}}
\put(2.2,3){\makebox(0,0)[ra]{$#7$}}
\put(14,5){\makebox(0,0)[la]{$#4$}} 
\end{picture}}
}   

\newcommand{\nbutterfly}[9]
{
{\Setunitlength
\begin{picture}(26,11)
\def\Kr{\circle{.4}}
\def\ID{\circle*{.4}}
\put(13,3){\Kr} 
\put(23,4){\ID}
\put(13,5){\ID}
\put(13,11){\Kr}
\put(11,8){\ID}
\put(15,8){\ID}
\put(5,7){\Kr}
\put(3,4){\ID}
\put(21,7){\Kr}
\put(11.1468,8.2202){\line(2,3){1.7065}}
\put(13.1468,5.2202){\line(2,3){1.7065}}
\put(14.8532,8.2202){\line(-2,3){1.7065}}
\put(12.8532,5.2202){\line(-2,3){1.7065}}
\put(3.1468,4.2202){\line(2,3){1.7065}}
\put(22.8532,4.2202){\line(-2,3){1.7065}}
\put(13.1788,3.0894){\line(2,1){7.6422}}
\put(3.1788,4.0894){\line(2,1){7.6422}}
\put(5.1788,7.0894){\line(2,1){7.6422}}
\put(12.8211,3.0894){\line(-2,1){7.6422}}
\put(22.8211,4.0894){\line(-2,1){7.6422}}
\put(20.8211,7.0894){\line(-2,1){7.6422}}
\put(23.1468,4.2202){\line(2,3){1}}
\put(23.1788,4.0894){\line(2,1){1}}
\put(2.8532,4.2202){\line(-2,3){1}}
\put(2.8211,4.0894){\line(-2,1){1}}
\put(1.5,6){\makebox(0,0)[rb]{$\cdots$}}
\put(24.5,6){\makebox(0,0)[lb]{$\cdots$}}
\put(13,11.8){\makebox(0,0)[cb]{$#1$}}
\put(13,2){\makebox(0,0)[ca]{$#9$}}
\put(15.8,8){\makebox(0,0)[lb]{$#3$}}
\put(10.2,8){\makebox(0,0)[rb]{$#2$}}
\put(21.8,7){\makebox(0,0)[lb]{$#6$}}
\put(4.2,7){\makebox(0,0)[rb]{$#5$}}
\put(23.8,3){\makebox(0,0)[la]{$#8$}}
\put(2.2,3){\makebox(0,0)[ra]{$#7$}}
\put(14,5){\makebox(0,0)[la]{$#4$}} 
\end{picture}}
}
 
\newcommand{\haekelmuster}[7]{
\noindent
{\Setunitlength 
\begin{picture}(15,9.5)
\def\Kr{\circle{.4}}
\def\ID{\circle{.4}}
\put(3,1){\Kr}
\put(3,4){\ID}
\put(3,7){\ID}
\put(6,1){\Kr}
\put(6,4){\ID}
\put(6,7){\ID}
\put(9,1){\Kr}
\put(9,4){\ID}
\put(9,7){\ID}
\put(12,1){\Kr}
\put(12,4){\ID}
\put(12,7){\ID}
\put(3,1.2){\line(0,2){2.6}}
\put(3,4.2){\line(0,2){2.6}}
\put(9,1.2){\line(0,2){2.6}}
\put(6,1.2){\line(0,2){2.6}}
\put(6,4.2){\line(0,2){2.6}}
\put(9,4.2){\line(0,2){2.6}}
\put(12,1.2){\line(0,2){2.6}}
\put(12,4.2){\line(0,2){2.6}}
\put(3.1414,3.8586){\line(1,-1){2.7172}}
\put(3.1414,6.8586){\line(1,-1){2.7172}}
\put(6.1414,3.8586){\line(1,-1){2.7172}}
\put(6.1414,6.8586){\line(1,-1){2.7172}}
\put(9.1414,3.8586){\line(1,-1){2.7172}}
\put(9.1414,6.8586){\line(1,-1){2.7172}}
\put(12.1414,3.8586){\line(1,-1){1}}
\put(12.1414,6.8586){\line(1,-1){1}}
\put(2.8586,1.1414){\line(-1,1){1}}
\put(2.8586,4.1414){\line(-1,1){1}}
\put(3,10){\makebox(0,0){$#7$}}
\put(1,4){\makebox(0,0)[rb]{$\cdots$}}
\put(14,4){\makebox(0,0)[lb]{$\cdots$}}
\put(3,8){\makebox(0,0)[cb]{$#1$}}
\put(3.4,4){\makebox(0,0)[lb]{$0$}}
\put(3.4,1){\makebox(0,0)[lb]{$0$}}
\put(6,8){\makebox(0,0)[cb]{$#2$}}
\put(6.4,4){\makebox(0,0)[lb]{$#4$}}
\put(6.4,1){\makebox(0,0)[lb]{$0$}}
\put(9,8){\makebox(0,0)[cb]{$#3$}}
\put(9.4,4){\makebox(0,0)[lb]{$#5$}}
\put(9.4,1){\makebox(0,0)[lb]{$#6$}}
\put(12,8){\makebox(0,0)[cb]{$0$}}
\put(12.4,4){\makebox(0,0)[lb]{$0$}}
\put(12.4,1){\makebox(0,0)[lb]{$0$}}
\end{picture}}
}
\newcommand{\posetDN}{
\noindent
{\Setunitlength 
\begin{picture}(15,8)
\def\Kr{\circle{.4}}
\def\ID{\circle{.4}}
\put(13,1){\Kr}
\put(13,4){\ID}
\put(13,7){\ID}
\put(1,1){\Kr}
\put(1,4){\ID}
\put(1,7){\ID}
\put(4,1){\Kr}
\put(4,4){\ID}
\put(4,7){\ID}
\put(7,1){\Kr}
\put(7,4){\ID}
\put(7,7){\ID}
\put(13,1.2){\line(0,2){2.6}}
\put(13,4.2){\line(0,2){2.6}}
\put(4,1.2){\line(0,2){2.6}}
\put(1,1.2){\line(0,2){2.6}}
\put(1,4.2){\line(0,2){2.6}}
\put(4,4.2){\line(0,2){2.6}}
\put(7,1.2){\line(0,2){2.6}}
\put(7,4.2){\line(0,2){2.6}}
\put(1.1414,3.8586){\line(1,-1){2.7172}}
\put(1.1414,6.8586){\line(1,-1){2.7172}}
\put(4.1414,3.8586){\line(1,-1){2.7172}}
\put(4.1414,6.8586){\line(1,-1){2.7172}}
\put(7.1414,3.8586){\line(1,-1){1}}
\put(7.1414,6.8586){\line(1,-1){1}}
\put(12.8586,1.1414){\line(-1,1){1}}
\put(12.8586,4.1414){\line(-1,1){1}}
\put(1.1940,1.0485){\line(4,1){11.6119}}
\put(1.1940,4.0485){\line(4,1){11.6119}}
\put(-3,7.7){\makebox(0,0)[lb]{$D_N:$}}
\put(9.5,5.5){\makebox(0,0)[lc]{$\cdots$}}
\put(9.5,2.5){\makebox(0,0)[lc]{$\cdots$}}
\put(13,7.7){\makebox(0,0)[lb]{$\zeta_{2p-3}$}}
\put(13.4,4){\makebox(0,0)[lb]{$\gamma_{2p-3}$}}
\put(13.4,.7){\makebox(0,0)[lb]{$\beta_{2p-3}$}}
\put(1,7.7){\makebox(0,0)[lb]{$\zeta_0$}}
\put(.5,4){\makebox(0,0)[rb]{$\gamma_0$}}
\put(.5,.7){\makebox(0,0)[rb]{$\beta_0$}}
\put(4,7.7){\makebox(0,0)[lb]{$\zeta_1$}}
\put(4.4,4){\makebox(0,0)[lb]{$\gamma_1$}}
\put(4.4,.7){\makebox(0,0)[lb]{$\beta_1$}}
\put(7,7.7){\makebox(0,0)[lb]{$\zeta_2$}}
\put(7.4,4){\makebox(0,0)[lb]{$\gamma_2$}}
\put(7.4,.7){\makebox(0,0)[lb]{$\beta_2$}}
\end{picture}}
}
\renewcommand {\epsilon}{\varepsilon}
\renewcommand {\leq} {\leqslant}
\renewcommand {\geq}{\geqslant}
\newcommand {\Oplus}{\bigoplus\limits}
\newdir { >} {{}*!/-5pt/@{>}}

\SelectTips{cm}{12} 
\title{Smash Products of $E(1)$-Local Spectra at an Odd Prime}
\author{Nora Ganter}
\date {\today}
\def\Setunitlength{\unitlength 12pt}
\begin{document}
\begin{abstract}
In \cite{Franke:95}, Franke constructed a purely algebraic
category that is equivalent 
as a triangulated category to the $E(n)$-local stable homotopy
category for $n^2+n<2p-2$.
The two categories are not Quillen equivalent, and 
his proof uses systems of triangulated diagram categories rather than
model categories.
Our main result is that in the case $n=1$ Franke's functor maps the
derived tensor product to the smash product. It can however not be
an associative equivalence of monoidal categories. 
The first part of our paper sets up a monoidal version of Franke's
systems of triangulated diagram categories and explores its properties.
The second part applies these results to the specific construction of
Franke's functor in order to prove the above result. 
\end{abstract}

\maketitle

\section{Background and Introduction}
\subsection{Chromatic Localizations}
The stable homotopy category is a very rich and complicated category,
and if homotopy theorists want to compute something, they often make
use of the fact that it has localizations that are easier to
understand. These localizations are typically Bousfield localizations
with respect to a homology theory. 
Localizations at the Johnson-Wilson homology theories 
$$E(n)_*(-)$$ 
have proven to be particularly fruitful for systematic
computations. These homology theories are defined 
by Landweber exactness of their coefficient groups
$$E(n)_* = \mathbf{Z}_{(p)}[v_1,\dots,v_n,v_n^{-1}]$$ 
over the Brown-Peterson spectrum $BP$. 
They are related to periodic phenomena in the stable homotopy groups
of spheres, therefore these localizations are often referred to as
chromatic localizations. 
The ``thick subcategory theorem" \cite{Hopkins:Smith} says that they
are all possible homology localizations of the category of finite
$p$-local spectra.
The spectrum $E(1)$ is also known as the Adams summand of the
$p$-local $K$-theory spectrum, and one has 
$$K_{(p)} = \bigvee_{i=0}^{p-2} \Sigma^{2i} E(1).$$ 
Therefore the localizations at $E(1)$ and at $K_{(p)}$ are the same.
\subsection{Franke's Algebraic Models}
From a categorical, structure theoretic point of view, the thick
subcategory theorem tells us that it will not be possible to find a model
for the $p$-local stable homotopy category that is as simple as
Serre's model for the rational stable homotopy category.
There is indeed a theorem by Schwede, asserting that
$\mathcal S_{(p)}$ has no ``exotic" model \cite{Schwede2}.
The localized categories $\mathcal S_{E(n)}$ promise to be simpler,
and indeed there are purely algebraic descriptions of these categories
by Bousfield and Franke.

Bousfield gives a purely algebraic description of the objects of the
$E(1)$-local stable homotopy category at an odd prime in
\cite{Bousfield:85}.
However, he says nothing about homotopy classes of maps.
He introduces an important tool, namely a specific setup for the Adams
spectral sequence by 
injective resolutions in the $E(1)$-local category: 
$$E_2^{s,t} = \Ext_{E(1)_*E(1)}^s(E(1)_*X,E(1)_*Y[t]) \Longrightarrow
[X,Y]_{E(1)\text{-local},t-s}.$$
Franke \cite{Franke:95}
uses Bousfield's work, in particular the computation of the injective
dimension of the category of $\EE$-comodules (it equals two) 
and the existence of this spectral sequence construction to 
show that $\mathcal S_{E(1)}$ is equivalent to the derived category of
a specific type of cochain complexes of $E(1)_*E(1)$-comodules.
He generalizes the result to ``higher chromatic primes", i.e.
$E(n)$-local spectra and $E(n)_*E(n)$-comodules (for $n^2+n\leq 2p-2$). 
\subsection{Systems of Triangulated Diagram Categories}
\label{tdcIntro}
It is important to note that Franke's functors are defined only on the
level of homotopy categories, and it can be shown that the
corresponding models are not Quillen equivalent.
This is the sense in which
Schwede refers to them as ``exotic" models.
The equivalences do however preserve the triangulated structure as well
as homotopy Kan extensions along maps of finite posets up to a certain
length. In particular they preserve homotopy (co)limits over such
diagrams. 
For the construction of Franke's functors it is also essential to have
functorial cones and homotopy Kan extensions at hand. How can one not
work on the level of model categories but still have a good handle on
homotopy Kan extensions? The answer is to work with homotopy
categories of strictly commuting diagrams of spectra; for each finite
poset one category. This is the idea of a system of triangulated
diagram categories. We explain it in more detail in Section
\ref{tdcausfuehrlich}.
In Section \ref{beweisvonthm} it becomes clear why this is the
natural setup for our project: Franke's functor is defined as a
homotopy colimit of a diagram of relatively simple shape. Once we
consider smash products, the diagrams become more complicated, and
many of our arguments are merely discussing the shape of the
underlying posets.
\subsection{Introduction}
Our paper falls into two major parts. First it describes in a general
setup how a monoidal structure, for example the smash product,
interacts with the structure maps of a system of triangulated diagram
categories. Then we apply our results to Franke's
functor in the case $n=1$. Denote this functor by $\Rek$.
The main result of the second part is 
\begin{Thm}\label{Maintheorem} There is a functor isomorphism
$$
\Rek(-) \wedge_{S_{E(1)}}^L \Rek(-) \cong
\Rek(- \otimes_{E(1)_*}^L -).
$$
\end{Thm}
Note that the category of $\EE$-comodules does not have enough
projectives and that therefore the existence of the derived tensor
product on the right hand side is a non-trivial statement.
In spite of our theorem, $\Rek$ is not a monoidal functor:
\begin{Rem}[Schwede]\label{Moore}
Let $p=3$, and denote the mod $3$ Moore spectrum by $M(3)$. 
This Moore spectrum has a unique multiplication which is not
associative. But
$$\Rek^{-1}(M(3)) \simeq 
\cdots\to0\to E(1)_*\stackrel{\cdot 3}{\to} E(1)_* \to 0\to\cdots $$
possesses an associative multiplication!
Thus we cannot hope to find an associative functor isomorphism between
$-\wedge^L-$ and the (derived) tensor product.
\end{Rem}
The proof of our theorem uses the specific construction of $\Rek$,
and many properties of systems of triangulated diagram categories.
Therefore we start by reviewing those parts of \cite{Franke:95}
that we need most frequently.
In order to make our work easily accessible for topologists, we stick
with the terminology of stable homotopy theory. However, 
we use only very general concepts, that could as well be formulated in
the language developed in \cite{Franke:95}.
\subsection{Plan}
In Section \ref{konstr} we recall the
parts from \cite{Franke:95} that we need for our computations. 
We discuss systems of triangulated diagram categories,
present the construction of the equivalence functor in 
a dual (but equivalent) version and state an easy generalization of
the spectral sequence \cite[1.4.35]{Franke:95}.
We do not assume familiarity with \cite{Franke:95}, but some of
our proofs use its results.  
Section \ref{diagrammsmash} is about the interaction of a monoidal
structure with the system of triangulated diagram categories.
It also contains some results about the computation of the edge
morphisms of a homotopy Kan extension, which turn out to be important
tools for our computations in Section \ref{beweisvonthm}.
The last part of Section \ref{diagrammsmash} can be skipped by a
reader who is only interested in the proof of Theorem
\ref{Maintheorem}. 
In Section \ref{welchesmodell} we explain which model for $\SE$ we choose
to work with. In Section \ref{flatmodel} we show that the derived
tensor product in the statement of Theorem \ref{Maintheorem} is well
defined and reduce the proof to the case of flat complexes.
Sections \ref{diagrammsmash}, \ref{welchesmodell} and \ref{flatmodel}
can be read after \ref{tdcausfuehrlich}. 
Finally, Section \ref{beweisvonthm} contains the proof of the main theorem.
\subsection{Acknowledgements} 
This paper is based on my Diplom thesis at the University of Bonn.
I would like to thank Jens Franke for
suggesting this thesis topic and for being an excellent advisor in every 
respect.  Many thanks also go to Stefan
Schwede and Birgit Richter for reading earlier drafts of this paper
and making plenty of helpful suggestions and to Haynes Miller, Mark
Hovey, Charles Rezk,
and Tilman Bauer for helpful comments and discussions. 
Last but not least I would like to take this opportunity to thank
Carl-Friedrich B\"odigheimer for helpful comments on this paper, but
also for being a wonderful teacher. His enthusiasm and his support for
young students meant a lot to me during my time in Bonn.

This research was partially supported by a Walter A.~Rosenblith 
fellowship and by a dissertation stipend from the German Academic 
Exchange Service (DAAD).
\section {Notations and conventions}
Let $p$ always be an odd prime. 
Whenever we draw a poset, the vertices represent elements,
and $x\leq y$ if and only if
the vertex corresponding to $x$ is linked to the vertex corresponding
to $y$ by an ascending path.
The length of a poset $C$ is defined to be
the supremum of all $k$ such that there exists a sequence
$x_0<x_1<\dots<x_k$ in $C$.
All posets we consider are finite and therefore of finite length.
Typically, if a cochain complex is called
$\C$, its cocycles are denoted $\Z$ and its coboundaries $\B$.
For our purposes it does not matter whether we work with symmetric spectra
\cite{HSS} or $S$-modules \cite{EKMM}. 
For a strict ring spectrum $R$, we denote the category of
strict $R$-module spectra by $\MM_R$ and its derived category by
$\Ho(\MM_R)$. More generally, $\MM$ always stands for a stable
model category \cite[7.1.]{Hovey}. 
If $E$ is a spectrum, we denote the underlying category of
$\MM_S$, endowed with the model structure used for Bousfield
localization 
at $E_*(-)$, by 
$$\MM_S \lbrack E^{-1}\rbrack.$$
Its derived category, i.e. the localization of the stable homotopy
category at $E_*(-)$ is denoted by
$$\mathcal S_E := \Ho(\MM_S\lbrack E^{-1}\rbrack).$$
Now
$$
\xymatrix{
{\MM_S}\ar @<.2ex> @{{}-^{>}} [r] ^{\id\phantom{iii}} & 
{\MM_S\lbrack E^{-1}\rbrack}\ar @<.2ex> @{{}-^{>}} [l] ^{\id\phantom{iii}}
}
$$
is a Quillen pair, and we write $(-)_E$ for the composition of its derived functors\footnote{Franke works with the localized categories, but topologists
  like to think of this functor as localization functor and refer to
  its image as $E$-local spectra.}
$$\mathcal S \stackrel{id^L}{\longrightarrow}
\mathcal S_E \stackrel{id^R}{\longrightarrow} \mathcal S.$$
It is induced by fibrant replacement in $\MM_S\lbrack E^{-1}\rbrack$
viewed as an endofunctor of $\MM_S$.
We will sometimes also write $(-)_E$ for the (strict) fibrant replacement
in $\MM_S\lbrack E^{-1}\rbrack$.
In particular, $S_E$ stands for the (strict)
$E$-local sphere, not as in \cite{EKMM} for the free $E$-module
spectrum (if $E$ is a ring spectrum). 
\begin {Def}
A quasi periodic cochain complex of period $N$ is a cochain complex 
$C_*^\bullet$ of
graded objects $C_*^n$ together with an isomorphism
$$C_{*-N}^\bullet = C_*^{\bullet}[N],$$
where the right hand side stands for $C_*^{\bullet}$ shifted to the
left $N$ times.
\end{Def} 
We denote the category
 of $\EE$-comodules which are concentrated
in degrees congruent to $0$ modulo $2p-2$ 
  by\footnote{In \cite{Bousfield:85} this category is called ${\mathcal
  B}(p)$, in \cite{Franke:95} it is $\tilde {\mathcal B}$.}
$$
\eeO.
$$ 
We denote the category of period $2p-2$ quasi periodic
cochain complexes in $\eeO$ by
$$\CN,$$ 
and its derived category 
by\footnote{Compare to Franke's notations $\mathfrak
C^{2p-2,[2p-2]}(\mathcal B)$ and 
  $\mathfrak D^{2p-2,[2p-2]}(\mathcal B)$.}
$$\mathcal D^{2p-2}(\eeO).$$
\section
{Franke's Algebraic Models}\label{konstr}
This section is a collection of the parts of \cite{Franke:95} needed
for our constructions.
We only consider the following special case of the main theorem of
\cite{Franke:95}:
\begin {Thm}[Bousfield, Franke]
The category of $E(1)$-local spectra is in a unique way equivalent to
the derived 
category of period $1$ quasi-periodic cochain complexes of
$E(1)_*E(1)$-comodules.
\end {Thm}
Here ``unique" means ``unique up to canonical natural isomorphism,
given that the equivalence is also valid for diagram categories for
diagrams of length $\leq 2$, preserves certain additional structure
between these, and transforms $E(1)_*(-)$ into something 
naturally isomorphic to $H^*(-)$''.
After we discuss the concept of a system of triangulated diagram
categories, we recall the construction of Franke's equivalence
functor. Finally we discuss a spectral sequence computing the homology
of homotopy Kan extensions. It is the major
computational device used throughout our paper.
\subsection{Systems of Triangulated Diagram Categories}
\label{tdcausfuehrlich}
The motivation for studying systems of triangulated diagram categories
was discussed in \ref{tdcIntro}. This section intends to give a more
detailed impression of this concept. We do not give the definition
here, but stick to the most important example:
The \elok is a triangulated category, and any of its standard 
models\footnote{See Section \ref{welchesmodell}.}
induces this triangulation. More precisely, a model category $\MM$ is
called {\it stable}, if the suspension functor is invertible in the
homotopy category $\Ho(\MM)$. This condition implies that $\Ho(\MM)$
is a triangulated category, where the triangles are the cofibre
sequences \cite[7.1]{Hovey}. 
Most triangulated categories in topology arise in this way.
Franke's paper never does any constructions on the level of model
categories. However,
the homotopy category alone is not
rigid enough. For example the cone is not a functor in a triangulated
category. Therefore, Franke also allows himself to work with homotopy
categories of categories of diagrams in $\MM$:
Let $C$ be a 
finite\footnote{One could do with much weaker conditions on $C$, 
  see \cite[5.2.]{Hovey}, \cite{DHK} and the original source
  \cite{Reedy}.} 
poset. The category $\MM^C$ of $C$ shaped diagrams in $\MM$ has a
model structure with vertex-wise weak equivalences and fibrations. The
cofibrations are characterized as 
follows\footnote{see e.g. \cite{Franke:95}. Here these cofibrations are
  called ``diagram cofibrations"}:
A morphism from $A$ to $B$ in $\MM^C$ is a Reedy cofibration if and
only if for all $c\in C$ the morphism
\begin{equation}\label{diagrammkof}
A_c\push{\colim{c'<c}A_{c'}}\colim{c'<c}B_{c'}\longrightarrow B_c,
\end{equation}
induced by the universal properties, is a cofibration in $\MM$.
Of course we could just as well define cofibrations and weak
equivalences vertex-wise and thus force the fibrations to be the maps
satisfying the dual condition to (\ref{diagrammkof}).
Note that both model structures have the same weak equivalences, and
therefore the same homotopy category. Franke always works on the
homotopy level, but he considers (for fixed $\MM$) the entire system
of homotopy categories of diagram model categories. This allows him to
use functorial homotopy Kan extensions, in particular
homotopy (co)limits, cones etc.: Let $f\negmedspace :C\to D$ 
be a map of posets. Our model structure is made in such a way that
pulling back diagrams along $f$ preserves fibrations and trivial
fibrations. Therefore the adjoint functors
$$
\xymatrix{
{\MM^C}\ar @<.2ex> @{{}-^{>}} [r] ^{\lkan f} & 
{\MM^D}\ar @<.2ex> @{{}-^{>}} [l] ^{f^*}
}
$$
form a Quillen pair, 
and $\lkan f$ has a left derived functor. Similarly, using the other
model structure, one defines right homotopy Kan extensions, right
adjoint to $f^*$.
The example of the cone functor \cite[1.4.5]{Franke:95} illustrates
how homotopy Kan extensions are computed and why they are useful.
\begin{Not}
Write $\II$ for the poset $\{ 0 < 1\}$, and
$$\VVv\subset\II\times\II$$
for the sub-poset with elements
$\{(1,0),(0,0),(0,1)\}$.
\end{Not}
%
%
\begin{Def}[Franke]\label{conedef}
Let $f\in\Ho(\MI)$. The cone of $f$ is defined as homotopy colimit
over the ``cokernel diagram" of $f$:
$$
{
\Setunitlength\raisebox{-1.75\unitlength}
{\begin{picture}(6,3.5)
\put(3,.5){\makebox(0,0)[ca]{$X$}}
\put(.5,3){\makebox(0,0)[cb]{$\star$}}
\put(5.5,3){\makebox(0,0)[cb]{$Y$}}
\put(4.5,1.4){\makebox(0,0)[lb]{$f$}}
\put(3.5,1){\line(1,1){1.5}}
\put(2.5,1){\line(-1,1){1.5}}
\end{picture}}
} 
$$
More precisely, let
$$\map{\Vst}{\MI}{\MM^\Vv}$$
denote the functor that takes as input an element $f$ of $\MI$ and returns
the object of $\MM^\Vv$ that has $f$ at the edge
$((0,0)\leq(0,1))$ and the (strict) zero object at the vertex $(1,0)$.
Then $\Vst$ preserves weak equivalences and thus induces a functor
$$\map{\Vst}{\Ho(\MI)}{\Ho(\Mv)}.$$
We define the cone of $f$ by
$$\cone(f) := \hocolim{\Vv}\Vst(f).$$
\end{Def}
\begin{Rem}
Note that $$\Vst \cong \horkan{\halfVv\subset\Vv}.$$
\end{Rem}
\begin{Rem}[Comparison with the classical cone definition]\label{kofkegel}
Let $\map fXY$ be a cofibration between cofibrant objects.
Let $CX$ be a cone object of $X$, i.e.
$$\xymatrix{X\ar@{>->}[r]&{CX}\ar@{->>}[r]^{\sim}&\star}.$$
We obtain a strict pushout diagram
\begin{equation}\label{classicalcone}
\xymatrix{
&CX\ar@{>->}[rd]\\
X\ar@{>->}[r]^f\ar@{>->}[ru] & Y \ar@{>->}[r] & {Y\push{f}CX.} 
}
\end{equation}
Without the bottom right corner, this is a Reedy cofibrant
replacement of $\Vst (f)$ in $\Mv$. Therefore the bottom right corner
is $\cone(f)$.
\end{Rem}
The functor $\cone$ takes values in $\Ho(\MM)$. Sometimes that is not enough.
Assume, for instance, we want to define cofibre sequences. Then we
have to take the cone of the cone inclusion map, i.e. of the bottom
right arrow in (\ref{classicalcone}). 
We need a functor
$$\map{\Cone}{\Ho(\MI)}{\Ho(\MI)}$$
that returns the bottom right arrow of (\ref{classicalcone}).
This is given by \cite[1.4.5]{Franke:95}
\begin{equation}\label{Conedef}
\Cone(-) := \left(\holkan{\Vv\subset(\Ii\times\Ii)}\circ\Vst(-)\right)
_{(0,1)<(1,1)}.
\end{equation}
\subsection{Construction of Franke's Functor}\label{frankesatz}
Using the functor $\E(-)$ Franke defines a ``reconstruction" functor,
$\Rek$, from the derived category of 
quasi periodic cochain complexes of $\EE$-comodules into the
category of $E(1)$-local spectra.
It turns out that $\Rek$ is an equivalence of categories. 
In order to motivate the way $\Rek$ is defined, we look at the
simplest analogous situation: instead of the system of diagram
categories corresponding to $E(1)$-local spectra,
we take the system corresponding to quasi periodic $\EE$-comodules 
and pretend that we want to reconstruct the identity functor by 
using only information that can be obtained via $H^*(-)$ (this now
plays the role of $\E(-)$).
Note that the data of a period $1$ quasi periodic cochain complex of
$\EE$-comodules are the same as the data of a period $2p-2$ quasi
periodic cochain complex of $\EE$-comodules that are concentrated in
degrees congruent to zero modulo $2p-2$, i.e.
\begin{equation}\label{periods}
{\mathcal C}^{1}(\operatorname{Comod}_{\EE})\cong\CN
\end{equation}
Now the basic idea is that any quasi periodic cochain complex $\C$ of
period $N$ can be decomposed into $N$ pieces of the form
\begin{equation}\label{Zn}
\xymatrix@=2ex{
\dots \ar[r] & 0 \ar[r] & C^n \ar[r] &
B^{n+1} \ar[r] & 0 \ar[r] & \dots \ar[r] & 0 \ar[r] & C^{n+N} \ar[r] &
B^{n+1+N}\ar[r] & 0 \ar[r] & \dots \\ \\
&& \save*[]  {\mbox{\tiny{$n^{th}$ spot}}}\restore
\ar@{..>}[-2,0],
}
\end{equation}
which can be glued back together along the inclusions of the $B^n$ into
$C^n$. 
More precisely, $\C$ is the colimit of the diagram of complexes
\begin{equation}\label{nKrone}
\xymatrix@=.1ex
{Z^0 &&&&&& Z^1           &&&&    && Z^2 &&&&&&&&&&&& Z^{N-1}\\
*+<4pt>[o][F]{} \ar@{-}[6,0] \ar@{-}[6,6] &&&&&& *+<4pt>[o][F]{} \ar@{-}[6,0]
\ar@{-}[6,6] &&&&&& *+<4pt>[o][F]{} \ar@{-}[6,0]\ar@{-}[3,3]
&&&&&&&&&&&& *+<4pt>[o][F]{} 
\ar@{-}[6,-24] \ar@{-}[6,0]\\  \\ \\                                   && &&
           &&&&&&&&                       &&&&   &&\dots&&&&&&\\ \\ \\
*+<4pt>[o][F]{} && && &&
*+<4pt>[o][F]{}                       &&&& && *+<4pt>[o][F]{} &&&&&&&&&&&&
*+<4pt>[o][F]{}\ar@{-}[-3,-3]\\ B^0 &&&&&& B^1
  &&&& && B^2 &&&&&&&&&&&& B^{N-1},\\
}\end{equation}
where we abbreviate (\ref{Zn}) by $Z^n$, 
and
$$\dots \to 0\to B^n \to 0 \to \cdots \to 0 \to B^{n+N}\to 0 \to\cdots$$ 
by $B^n$;
the vertical edges are $B^n\hookrightarrow C^n$, 
and the diagonal edges are $B^{n+1} = B^{n+1}$. 
Note that the colimit of the diagram (\ref{nKrone}) is equal to its
homotopy colimit by (\ref{diagrammkof}).  
How can we read off $\C$ from the cohomology of such a diagram?
Note first that
$$\dots \to 0 \to C^n \to 0 \to \dots$$
turns up as cone of the diagonal maps.
The cone inclusion to the complex $\Sigma B^{n+1}$ is $d^n$.
(Note that $B^{n+1}$ is concentrated in degrees congruent
to $n+1$, and that $Z^n$ and $C^n$ are concentrated in degrees
congruent to $n$.)
But $H^*$ reflects
$$\xymatrix {
\dots \ar[r] & 0 \ar[d]\ar[r] & C^n \ar@{->>}[d]^d \ar[r] &
0 \ar[d]\ar[r] & \dots &&  C^n\ar[d]\\
\dots \ar[r] & 0 \ar[r] & B^{n+1} \ar[r] &
0 \ar[r] & \dots && \Sigma B^{n+1}.\\
}$$
And the same is true for the composite 
$$\xymatrix {
\dots \ar[r] & 0 \ar[d]\ar[r] & B^{n+1}
\ar@{>->}[d]\ar[r] & 0\ar[d]\ar[r] & \dots &&\\
\dots \ar[r] & 0 \ar[d]\ar[r]
& C^{n+1} \ar@{=}[d]\ar[r] & B^{n+2}\ar[d]\ar[r] &0\ar[r]& \dots \\
\dots
\ar[r] & 0 \ar[r] & C^{n+1} \ar[r] & 0 \ar[r] & \dots &&\\
}
\xymatrix{
B^{n+1}\ar[d]^{\parbox{30mm}{\raggedright\tiny vertical edge}}\\
Z^{n+1}\ar[d]^{\text{cone inclusion}}\\
C^{n+1}
}$$
In other words, it is possible to reconstruct $C^\bullet$ from
(\ref{nKrone}) by applying $H^*$. 
We return to $E(1)$-local spectra.
\begin{Not}
We denote the underlying poset of (\ref{nKrone}) by $C_N$. It has vertices
$\zeta_n$ and $\beta_n$, where 
$n\in\mathbb{Z}_{\!/\!\raisebox{-.2ex}{$_N$}}$ and relations
$\beta_n \leq \zeta_n$ and $\beta_{n+1} \leq \zeta_n$. 
In our case, $N = 2p-2$. 
\end{Not}
We consider the full subcategory
of objects $A$ of $\Ho(\EM^{C_N})$ satisfying:
\begin {itemize}{\sloppy
\item $Z^n_* := E(1)_{*-n}(A_{\zeta_n})$ and $B^n_* := E(1)_{*-n}(A_{\beta_n})$
are concentrated in degrees \linebreak[3]\mbox{$\equiv 0 \hbox{ $\mod 2p-2$}$}, and
\item $\map {E(1)_{*-n}(A_{\beta_n} \to A_{\zeta_n})} {B^n_*} {Z^n_*}$ is
injective.
}\end {itemize}
\begin{Not}\label{LNot}
Following Franke, we denote this subcategory of $\Ho(\EM^{C_N})$ by
$\mathcal L$.
\end{Not}
\begin{Con}\label{QConst}
Let $A$ be an object of $\mathcal L$. We define
$$C^n_*(A) := E(1)_{*-n}\left(\cone(A_{\beta_{n+1}}
\to A_{\zeta_n})\right).$$ 
If we apply $E(1)_{*-n}(-)$ to the exact triangle
$$
A_{\beta_{n+1}} \to A_{\zeta_n} \to \cone \left(A_{\beta_{n+1}}
\to A_{\zeta_n}\right) \to\Sigma A_{\beta_{n+1}} \to\Sigma A_{\zeta_n}
$$
we obtain a (short) exact sequence
$$B^{n+1}_{*+1} \overset 0 \rightarrow Z^n_* \rightarrow C^n_* 
\rightarrow B^{n+1}_*
\overset 0 \rightarrow Z^n_{*-1}.$$
It follows that $C^n_*$ is also concentrated in degrees $\equiv 0
\hbox{ $\mod 2p-2$}$. 
In order to define the differential, we apply $E(1)_{*-n}(-)$ to
$$\cone \left(A_{\beta_{n+1}} \to A_{\zeta_n}\right) \to
\Sigma A_{\beta_{n+1}} \to A_{\zeta_{n+1}} \to \cone\left(A_{\beta_{n+2}} 
\to A_{\zeta_{n+1}}\right)$$
and obtain
\begin{equation}\label{diff}
\xymatrix {
C^n_*(A) \ar@{->>}[r]\ar@/_4ex/[0,3]_{d^n.} & 
B^{n+1}_*\phantom{|} \ar@{>->}[r] & 
Z^{n+1}_*\phantom{|} \ar@{>->}[r] & C^{n+1}_*(A)}
\end{equation}
\end{Con}
Franke proves that this defines an equivalence of categories
\begin{eqnarray*}
\mathcal L & \to & \CN\\
A &\mapsto & C_*^\bullet (A).
\end{eqnarray*} 
\begin{Not}\label{QNot}
We call this equivalence $Q$.
\end{Not}
Let
$Q^{-1}$ be an inverse of $Q$. 
Franke further shows that 
\begin{equation}\label{Recdef}
\map {\Rek := \hocolim{} \circ Q^{-1}}{\CN}{\SE}
\end{equation}
factors over the derived category 
and induces an equivalence of categories 
$$\DD^1(\text{$E(1)_*E(1)$-comod}) \simeq \DD^{2p-2}(\eeO ) \to \SE,$$
the first equivalence being (\ref{periods}).
We denote this equivalence also by $\Rek$.
\begin{Rem}
Franke's construction uses cokernels and coimages, rather than kernels
and images. The reason for this is that we need to use Adams spectral
sequences via injective resolutions --- there are not enough projectives.
The theorem is also valid for 
diagram categories up to a certain length of the diagrams. In
order to improve the restrictions on this length, it is advantageous
to work with cokernels and coimages. However, Franke also proves a
uniqueness statement about $\Rek$, which implies that the dual construction
we have discussed here gives us the same functor on $\SE$.
\end{Rem}
\subsection{A spectral sequence}\label{spektralfolge}
There is yet another construction we will use from Franke's article,
namely the spectral sequence \cite[1.4.35]{Franke:95}. 
We rephrase his proof to give a slight generalization for homotopy Kan
extensions: 
Let $\MM$ be a stable model category, let $f\negmedspace :D\to C$ 
be a map of finite posets, and let $Y\in\Ho(\MM^D)$. 
Further, let 
$$F_*\negmedspace :\Ho(\MM)\to\mathcal A$$ 
be a homological functor into a Grothendieck category of graded objects,
such that
$$F_{t+1}(X) = F_t(\Sigma X).$$ 
Then there is a spectral sequence
\begin{equation}\label{SF}
E_2^{s,t}={\lkan f}{}_{-s}F_t(Y)\Longrightarrow
F_{s+t}(\holkan{f} Y).
\end{equation}
The example that is relevant for us is
$$\map{F_t(-) := E(1)_{-t}(-)}{Ho(\EM)}{\eee}.$$ 

The {\bf Construction}
of (\ref{SF}) is analogous to that of \cite[1.4.35]{Franke:95}:
Let $D$ be a poset. For $d\in D$, let $\map{i_d}{\{ d\}}{D}$ denote the 
inclusion of the vertex $d$ in $D$.
For $Y\in\Ho(\MM^D)$, define
$$\PP Y := \bigoplus_{d\in D}\holkan{i_d} Y_d .$$ 
Then we have:
\begin{equation}
  \label{tildeP-eqn}
  (\PP Y)_{d'} = \Oplus_{\{d\in D \mid d\leq d'\}}Y_d. 
\end{equation}
We define a morphism $g\negmedspace :\tilde PY\to Y$ by letting
$$\map{g\arrowvert_{\holkan{i_d}Y_d}}{\holkan{i_d}Y_d}{Y}$$ 
be the counit of the adjunction. Taking the homotopy fibre $\tilde RY$
of $g$
and then iterating the whole process with $\tilde RY$ in the role of $Y$,
we obtain a resolution
\begin{equation}\label{auflosung}
\xymatrix{
\cdots\ar@{=}[r]&{\Sigma Y}\ar@{=}[r]^{+}&Y\ar[r]^{+}\ar[ld]&\tilde RY
\ar[r]^{+}\ar[ld]&\tilde R^2Y\ar[r]^{+}\ar[ld]&\dots\ar[r]^{+}&\tilde
R^{\hgt(D)}Y\ar[r]^{+}&0\ar[ld]\\
               & 0\ar[u]        &\PP
Y\ar[u]^g    &\PP\tilde RY\ar[u]            &  &  &{\PP\tilde R^{\hgt(D)}Y}
\ar@{=}[u]   &{\phantom{,}0.}\ar@{=}[u]_{\parbox{1cm}{ $\cdots$}}\\ 
}
\end{equation}
Since
$F_*(\holkan{f}(-))$ is a homological functor, we obtain an
exact couple 
$$\xymatrix{D_1^{**}\ar[0,2]^{(-1,1)}_\alpha & &
D_1^{**}\ar[2,-1]^{(1,0)}_\beta\\
\\
&E_1^{**}\ar[-2,-1]^{(0,0)}_\gamma
\\
}$$
where
\begin{eqnarray*}
D_1^{s,t}&=&F_{t}(\holkan{f} (\tilde R^{-s}Y))\\
E_1^{s,t}&=&F_{t}(\holkan{f} (\PP\tilde R^{-s}Y)).
\end{eqnarray*}
This gives rise to a (cohomological) spectral sequence in the usual way.

We want to discuss its {\bf Convergence}: Let
$r>\hgt D$. After we derive our exact couple $r-1$ times, the
resulting couple
$(D_r^{**},E_r^{**})$ is:
$$\begin{array}{l}
D_r^{s,t}=\im(\map{\alpha^{r-1}}{D_0^{s+r-1,t-r+1}}{D_0^{s,t}}) \\
\\
=\left\{
\begin{array}{ll} D_0^{0,s+t}         & \mbox{if } s\geq 0,     \\
\im(\alpha^{-s})     & \mbox{if } -\hgt D\leq s<0, \\0 &
\mbox{else.}
\end{array}\right\}
\end{array}$$
It is of the form
$$\xymatrix{
\cdots D_0^{1,*+1}\ar@{=}[r]^{+}&D_0^{0,*}\ar@{->>}[r]^{+}\ar[ld]^{0}&
  \im (\alpha)\ar@{->>}[r]^{+}\ar[ld]^{0}
&\im(\alpha^2)\ar@{->>}[r]^{+}\ar[ld]^{0}&\dots\ar@{->>}[r]^{+}&
\im(\alpha^{\hgt(D)}) \ar@{->>}[r]^{+}&0\ar[ld]^0\cdots\\
0\ar@{>->}[u] &E_r^{0,*}\ar@{>->}[u]   &E_r^{-1,*}\ar@{>->}[u]
&E_r^{-2,*}\ar@{>->}[u]&&E_r^{-\hgt(D),*}.\ar@{>->}[u] \\
}$$
Therefore the spectral sequence collapses after $(\hgt (D)+1)$ steps,
and it converges as follows:
$$\xymatrix@=2ex{
F_{t}(\holkan{f}
Y)\ar@{->>}[0,2]&&D_\infty^{-1,t+1}\ar@{->>}[0,2]&
&D_\infty^{-2,t+2}\cdots\ar@{->>}[0,1]&
D_\infty^{-\hgt D,t+\hgt D}\ar@{->>}[0,2]&&0\\
\\
&\save *[]{\ker = E_\infty^{0,t}}\restore\ar@{..>}[-2,0] &&\save *[]
{\ker = E_\infty^{-1,t+1}}\restore\ar@{..>}[-2,0]&&\dots\\}$$

\noindent For the {\bf{identification of the $E_2$ term}}
we need to show that
\begin{equation}\label{resolution-eqn}
  F_{t}(Y)\twoheadleftarrow
  F_{t}(\PP\tilde R^\bullet Y)
\end{equation} 
is a $\lkan{f}$-acyclic resolution.
It is clearly exact. 
That it is a resolution by $\lkan{f}$-acyclic objects follows from the
equality
$$F_{t}(\PP Y') = \Oplus_{d\in
D}\underset{i_d}{\operatorname{LKan}}F_{t}(Y'_d)$$ 
together with the following lemma.
\begin{Lem}\label{azyklisch}
Let ${\mathcal A}$ be a Grothendieck category. 
Let $D$ be a finite poset
and assume that $X\in {\mathcal A}^D$ is such that for any $d\in D$ 
the map
$$
\colim{c< d}X_c\rightarrow X_d,
$$
given by the universal property of the colimit applied to the
edges of $X$, is a monomorphism. Then for
any map of finite posets $f\negmedspace : D\to C$, the object $X$ is 
$\lkan{f}$-acyclic. 
In particular, $X$ is also $\colim{}$-acyclic.
\end{Lem}
\begin{Pf}{}
We endow the category $\mathcal C(\mathcal A)$ 
of cochain complexes in $\mathcal A$ 
with the injective model 
structure\footnote{By that we mean that cofibrations are degree-wise
  monomorphisms, weak equivalences are the quasi isomorphisms, and
  fibrations are degree-wise epimorphisms with degree-wise injective
  cokernel. For the existence of such a model structure on the
  category of unbounded cochain complexes in a  
  Grothendieck-category, see \cite{Franke:Brown} or  \cite{Hovey:Groth}.
  There one can also find a discussion of the connection between this
  story and the derived functors in the sense of homological algebra.
}.
For the category of $D$-diagrams of chain complexes,
$$\mathcal C(\mathcal A^D)\cong \mathcal C(\mathcal A)^D,$$
we choose the model structure, whose fibrations and weak equivalences
are defined vertex-wise, and whose cofibrations are characterized by
(\ref{diagrammkof}). 
Then the satellite functors of the total derived functor
$$\mathcal A^D \to \mathcal D(\mathcal A^D)
\stackrel{\holkan f}{\longrightarrow}\mathcal D(\mathcal A^C)$$
satisfy the universal property of the (left-) derived left Kan
extensions. Now the condition of the lemma on $X$ says exactly that
$X$, viewed as an object of $\mathcal C(\mathcal A^D)$ (i.e.\ as a
cochain complex which is
concentrated in degree zero), is Reedy cofibrant. 
Thus in $\mathcal D(\mathcal A^D)$, we have
$$\holkan fX\cong \lkan fX.$$
The object $\lkan fX$, however, is also concentrated in degree
zero. Therefore 
all higher derived functors vanish.
\end{Pf}
We have shown that the objects in the resolution (\ref{resolution-eqn})
are $\lkan f$-acyclic. 
Further, by (\ref{tildeP-eqn}), we have
$$F_{t}(\holkan{f}\PP Y')_c =
F_{t}(\Oplus_{f(d)\leq c} Y'_d) = \Oplus_{f(d)\leq c
}F_{t}(Y'_d) =\left(\lkan{f} F_{t}(\PP Y')\right)_c.$$
Therefore we have shown that the cohomology of the complex
$F_{t}(\holkan{f}\PP\tilde R^\bullet Y)$ 
computes the derived functors of the left Kan extensions.

In calculations, we sometimes write down the $E_2$-term
directly.
If the reader is not familiar with computing derived Kan extensions, (s)he can
check their correctness by writing down the $E_1$-term.
\section{Smash products for diagram categories}\label{diagrammsmash}
This section discusses the interaction of a monoidal structure with a
system of triangulated diagram categories.
Our first goal is to define a smash product between the (homotopy)
diagram categories.
We start by defining a strict smash product of diagrams.
\begin{Def}Let $(\MM,-\wedge-)$ be a (model) category with monoidal structure.
Let $C$ and $D$ be finite posets, let
$$f\negmedspace : X\to Y \in\Mor(\MM^C)\text{  and  } 
g \negmedspace : U\to V \in\Mor(\MM^D).$$
We define $X \wedge U \in \Ob(\MM^{C \times D})$ by 
$$(X \wedge U)_{(a,b) \leq (c,d)} := X_{a \leq c} \wedge U_{b \leq d}$$ 
and $f\wedge g\in\Hom_{\MM^{C \times D}}(X\wedge U,Y\wedge V)$ by
$$(f\wedge g)_{(c,d)} := f_c\wedge g_d.$$
\end{Def}
In the following we ask for the monoidal structure to be
compatible with the model structure and conclude that in this case also
the smash products of diagrams are compatible with the model
structure. 
By ``compatible" we mean in the sense of the definitions 
\cite[4.2.]{Hovey}. 
For our diagram categories these definitions read:
\begin{Def}
In the situation of the previous definition
  the pushout smash product of $f$ and $g$ is defined to be the
  canonical map
  $$\map{f \Box g}{X\wedge V\push{X\wedge U} Y\wedge U}{Y\wedge V}
  \in\Mor(\MM^{C\times D}).$$
\end{Def}
\begin{Rem}
  If we view $f$ and $g$ as objects of $\MM^{C\times\Ii}$ and 
  $\MM^{D\times\Ii}$,
  we have
  \begin{equation}\label{box}
    {f\Box g = \lkan{\id_{C\times D}\times\pC} f\wedge g},
  \end{equation}
  where $\pC$ is the map
  \begin{eqnarray*}
    \pC\negmedspace : \II\times\II & \to     &   \II \\
    (1,1)                          & \mapsto &  1    \\
    (1,0),(0,0),(0,1)              & \mapsto &  0.
  \end{eqnarray*}
\end{Rem}
\begin{Def}
We say that
$$\map{-\wedge -}{\MM^C\times\MM^D}{\MM^{C\times D}}$$
satisfies the monoid axiom, if for any two cofibrations $f$ in $\MM^C$
and $g$ in $\MM^D$, the map $f\Box g$ is a cofibration, which is
trivial if $f$ or $g$ is.
\end{Def}
\begin{Def}
A (symmetric) monoidal model category is a model category $\MM$
together with a closed (symmetric) monoidal structure 
$-\wedge -$, such that
\begin{itemize}
\item there are functorial factorizations of morphisms into a trivial
cofibration followed by a fibration, and into a cofibration followed
by a trivial fibration respectively,
\item $\map{-\wedge -}{\MM\times\MM}{\MM}$ satisfies the monoid axiom,
and
\item the cofibrant replacement of the unit
$$\xymatrix{{QS}\ar[r]^{q}&S}$$ satisfies: for any cofibrant object $X$,
$$\xymatrix{QS\wedge X\ar[r]^{q\wedge\id_X}&S\wedge X}$$
is a weak equivalence (and
$$\xymatrix{X\wedge QS\ar[r]^{\id_X\wedge q}&X\wedge S}$$
is, too).
\end{itemize} 
\end{Def}
Hovey showed in \cite[4.3.1.]{Hovey} that for such a (symmetric)
monoidal model category, the monoidal structure on $\MM$ has a
left derived functor, which is itself a (symmetric) monoidal structure
on the homotopy category.
For our application, only one of the following two examples will be relevant
(compare Section \ref{welchesmodell}).
\begin{Exa}\label{Rmod}{}
  The model category $\MM_R$ of strict modules over a strict, strictly
  commutative ring spectrum $R$ is a symmetric monoidal model category
  \cite{SS:97}. 
\end{Exa}
\begin{Exa}\label{BL}{(I learned this from Stefan Schwede.)}
  Let $(\MM_S,-\wedge_S-)$ be a model for the stable homotopy
  category, and let $E$ be an object of $\MM_S$. Then
  $$(\MM\lbrack E^{-1}\rbrack,-\wedge_S-)$$ 
  is also a monoidal model category.
\end{Exa}
\begin{Pf}{}
  Without loss of generality we may choose $E$ to be cofibrant.
  Since Bousfield localization does not change the cofibrations,
  we only have to check the monoid axiom for
  $$
    f\negmedspace :X \to Y \text{ and  } 
    g\negmedspace :U \to  V 
  $$
  in the case that $g$ is an $E$-isomorphism, i.e. if
  $$ g\wedge\id_E:\negmedspace U\wedge E \to V\wedge E $$
  is a week equivalence.
  The monoid axiom for $(\MM_S,-\wedge_S-)$ implies that
  $$g\wedge\id_E = g\Box
  (\xymatrix{\star\ar@{>->}[r] & E})$$ 
  is a cofibration, because $E$ is cofibrant.
  If we apply the monoid axiom for $$(\MM_S,-\wedge_S-)$$ once more,
  this time to $f$ and $g\wedge\id_E$, it follows that
  $$\map{f\Box (g\wedge\id_E) = (f\Box g)\wedge\id_E}
    {\left(X\wedge V\push{X\wedge U}Y\wedge U\right)
    \wedge E}{Y\wedge V\wedge E}$$
  is also a weak equivalence.
\end{Pf}
\begin{Prop}\label{knall}{}
Let $(\MM,-\wedge-)$ be a monoidal model category. Then
$$\map {- \wedge-} {\MM^C \times\MM^D} {\MM^{C\times D}}$$ 
has a total left derived functor
$$\map{- \wedge^L -} {\Ho(\MM^C) \times \Ho(\MM^D)} {\Ho(\MM^{C\times D})}.$$
\end{Prop}
\begin{Pf}{}
By \cite[4.3.1.]{Hovey}, it is enough to show that 
the monoid axiom is satisfied:
Let $U \rightarrowtail V$ be a Reedy cofibration in $\MM^C$
and $X\rightarrowtail Y$ a Reedy cofibration in $\MM^D$. 
I.e. we assume
\begin{equation}\label{kofibrantinM}
\forall c \in C: U_c \push{\colim{a < c} U_a} \colim{a<c} V_a
\rightarrowtail V_c 
\end{equation}
and
$$\forall d \in D: X_d \push{\colim{b <
d} X_b} \colim{b<d} Y_b \rightarrowtail Y_d$$
to be cofibrations in $\MM$. Since
$- \wedge A$ and $A \wedge -$ are left adjoints and therefore preserve
colimits, and because of the monoid axiom in
$\MM$, it follows from (\ref{kofibrantinM}) that
the morphism with source the pushout of
$$U_c \wedge Y_d \push{\colim{a<c}U_a \wedge
Y_d}  \colim{a<c}V_a \wedge Y_d
\quad\text{and}\quad 
V_c\wedge  X_d\push{\colim{b<d}V_c \wedge X_b} \colim{b<d}V_c \wedge Y_b$$
over
$$\left(U_c  \push{\colim{a<c}U_a}\colim{a<c}
V_a\right) \wedge \left(X_d  \push{\colim{b<d}X_b} \colim{b<d} Y_b\right)$$
and target $V_c \wedge Y_d$ is a cofibration for every pair $(c,d)$.
By the remark below this pushout is just the pushout of
$$U_c \wedge Y_d \push{U_c\wedge X_d} V_c\wedge X_d \quad\text{and}\quad
\colim{(a,b)<(c,d)}V_a \wedge Y_b $$ over 
$$\left(\colim{(a,b)<(c,d)} U_a \wedge Y_b\push{\colim{(a,b)<(c,d)} U_a
\wedge X_b} \colim{(a,b)<(c,d)} V_a \wedge X_b\right),$$ 
and the map is the one you would expect there.
We have shown that
$$U\wedge Y\push{U\wedge X}V\wedge Y\to V\wedge Y$$ 
is a Reedy cofibration.
If one of the two maps
$U \rightarrowtail V$ and $X \rightarrowtail Y$ 
is a vertex-wise weak equivalence, it follows from the monoid axiom in
$\MM$ that
$U\wedge Y\push{U\wedge X}V\wedge Y\to V\wedge Y$ 
is also a vertex-wise weak equivalence.
Therefore the monoid axiom is satisfied.
\end{Pf}
\begin{Rem}
The calculation uses the isomorphism
$$\colim{(a,b)<(c,d)}?_a \wedge ?_b \cong \colim{a<c} ?_a \wedge
?_d\push{\colim{\underset{b<d}{a<c}}?_a \wedge ?_b} \colim{b<d}?_c \wedge
?_b,$$ 
which follows in a straightforward way from the various universal
properties. 
Alternatively, by (\ref{striktekanten})
the right hand side is
$$
\colim\Vv\lkan{f} (?_a\wedge ?_b)
$$ 
with
\begin{eqnarray*}
f \negmedspace : \{(a,b)\mid(a,b)<(c,d)\} & \to & \VVv \\
(c,b) & \mapsto & (1,0) \\
(a,d) & \mapsto & (0,1) \\
(a,b) \mid a<c \text{ and }b<d &\mapsto & (0,0) 
\end{eqnarray*}
(compare notation 1).
\end{Rem}
\begin{Rem}
In the proof of the preceding proposition we are working with the
model structure that has Reedy cofibrations as cofibrations.
We do so for later reference.
Of course, if we work with vertex-wise cofibrations and weak
equivalences instead, the monoid axiom is straight forward.
The universal property of $-\wedge^L -$
implies that up to canonical isomorphism both constructions give
the same result.
\end{Rem}
In the following, $(\MM,-\wedge-)$ is a monoidal,
stable model category.
We discuss compatibility of $-\wedge^L-$ with the various
other structures of the system $\Ho(\MM^C)$.
In order to do so, we need a lemma about the composition of derived
functors.
\begin{Lem}\label{GoF}
Let $\mathcal C$ and $\mathcal D$ be model categories, $\mathcal E$ an
arbitrary category, and let
$$\map{F}{\mathcal C}{\mathcal D} \quad\text{and} \quad 
\map{G}{\mathcal D}{\mathcal E}$$
be functors, such that $F$ sends trivial cofibrations between
cofibrant objects to weak equivalences between cofibrant
objects,
and $G$ sends trivial cofibrations between cofibrant objects to isomorphisms.
Then the derived functors $LG$, $LF$ and $L(G\circ F)$ exist, and 
in the commutative diagram of natural transformations
\begin{equation}\label{LGoF}
\xymatrix{
{LG\circ LF} \ar[0,2]^{\iota}\ar[1,1] & & {L(G\circ F)}\ar[1,-1] \\
& G\circ F ,&\\
}
\end{equation}
induced by the universal properties of the various derived functors,
$\iota$ is a functor isomorphism. Moreover, $\iota$ is associative up
to canonical natural equivalence.
\end{Lem}
\begin{Pf}{}
The existence of the derived functors is for example discussed in
\cite{Hovey}.
By construction, the diagonal morphisms in
(\ref{LGoF}) are isomorphisms on cofibrant objects.
Therefore $\iota_X$ is an isomorphism for cofibrant $X$.
For arbitrary $X$ we conclude
$$
\xymatrix{
LF\circ LG(X)\ar[r]^{\iota_X} & L(G\circ F)(X) \\
LF\circ LG(QX)\ar[r]^{\iota_{QX}}_{\cong}\ar[u]_{\cong} & L(G\circ F)(QX)
\ar[u]_{\cong}, \\
}
$$
where for the moment $Q(-)$ denotes the cofibrant replacement functor.
The argument also shows the associativity of $\iota$.
\end{Pf}
\begin{Exa}
Left Quillen functors preserve trivial cofibrations, cofibrations and
initial objects. Therefore they also preserve cofibrant objects (and
trivial cofibrations between them).
\end{Exa}
\begin{Cor}\label{smhol}
There is a functor isomorphism
$$\hocolim {C \times D}
(A \wedge^L B) \cong (\hocolim C A) \wedge^L (\hocolim D
B).$$
\end{Cor}
\begin{Pf}{}
As in \cite[4.3.1]{Hovey},
it follows from the monoid axiom (proof of Proposition \ref{knall})
that $-\wedge -$ 
sends Reedy cofibrant objects of $\MM^C\times\MM^D$ to diagram
cofibrant objects of $\MM^{C\times D}$ and preserves
(vertex-wise) trivial Reedy cofibrations between Reedy cofibrant
objects.
As left a Quillen functor, $\colim{}$ also
satisfies the conditions of the lemma on $F$. 
Further, for $A\in\MM$, we know that $A\wedge -$ is a left adjoint and
therefore commutes with colimits, which implies the analogous strict
formula.  
\end{Pf}
More generally, we have
\begin{Cor}
There is a functor isomorphism
$$
\holkan{f\times g} (A\wedge^L B)\cong (\holkan f A)\wedge^L(\holkan gB). $$
\end{Cor}
\begin{Pf}{}
The proof is analogous to the proof of Corollary \ref{smhol}. The
strict formula
follows, because left Kan extensions commute with left adjoints.
\end{Pf}
\begin{Cor}\label{smzur}
There is a functor isomorphism
$$f^* A\wedge^L g^* B \cong (f\times g)^*(A\wedge^LB).$$
\end{Cor}
\begin{Pf}{}
This time we work with the model structure whose cofibrations and weak
equivalences are defined vertex-wise.
Then $-\wedge-$ and pulling back both preserve cofibrant objects and
trivial cofibrations between them.
The analogous strict statement
follows directly from the definition. 
\end{Pf}
\begin{Cor}
The pushout smash product has a left derived functor
$$\map{-\Box^L-}{\Ho(\MI)\times\Ho(\MI)}{\Ho(\MI)}.$$
\end{Cor}
\begin{Pf}{}
According to (\ref{box}), we have
$$f\Box g = \lkan{\pC}(f\wedge g).$$
But
$$\map{-\wedge-}{\MI\times\MI}{\MM^{\Ii\times\Ii}}$$
and $\lkan{\pC}$ both preserve Reedy cofibrant objects and
vertex-wise trivial Reedy cofibrations between them.
Therefore the left derived pushout smash product exists and is given
by 
\begin{equation}\label{lbox}
-\Box^L- = \holkan{\pC}(-\wedge^L-).
\end{equation}
\end{Pf}
There is one further corollary, that has nothing to do with the
monoidal structure, but will turn out to be useful.
\begin{Not}\label{KantenNot}
Let $\map{f}{C}{D}$ be a map of posets.
For $d\in D$, we let $C\to d$ denote the subposet
$$\{c\in C\mid f(c)\leq d\}$$ 
of $C$, and we let 
$$\map{j_d}{(C\to d)}{C}$$ 
denote its inclusion into $C$.
For an edge $d\leq d'\in D$, let
\begin{eqnarray*}
p_d^{d'}\negmedspace : (C\to d')& \longrightarrow & \II \\ 
c &\mapsto & 0 \quad\text{if $f(c)\leq d$,} \\
c & \mapsto & 1 \quad\text{else.}
\end{eqnarray*}
\end{Not}
The vertices of a left
homotopy Kan extension are given by \cite[Prop.1.4.2]{Franke:95}:
\begin{equation}\label{ecken}
(\holkan{f}(X))_d \cong 
\hocolim{C\to d}j_d^*X.
\end{equation}
In an algebraic situation, where it makes sense to speak about
satellite functors, this becomes
\begin{equation}\label{striktekanten}
\left(\lkan f{}_sX\right)_d \cong \left(\colim{C\to d}\right)_s
\left(X\at{C\to d}\right).
\end{equation}
The following corollary says that the edges of a homotopy Kan
extension are also what we would expect. 
\begin{Cor}[Edges of homotopy Kan extensions]\label{Kanten}
There is a functor isomorphism
$$(d\leq d')^*\holkan fX\cong\holkan{p_d^{d'}}j_{d'}^*X.$$
\end{Cor}
This implies that the edges of 
$$\left(\lkan f{}_sX\right)$$
are given by the various universal properties on the right hand side of
(\ref{striktekanten}).

\begin{Pf}{}
Look at the analogous strict diagram, its vertices are given by
(\ref{striktekanten}) with $s=0$,
and at the edge
$d\leq d'$ there is the map that is obtained by applying the universal
property of
$\colimlim{c\in C\to d}X_c$
to the colimit inclusions
$$X_c\hookrightarrow\colim{c'\in C\to d'}X_{c'}.$$  
This diagram satisfies the universal property of the strict left Kan
extension of $X$ along $f$ in $\MM^C$. 
Therefore the analogous strict statement is true.
Since  $j_{d'}^*$ and left Kan extensions preserve diagram
cofibrations, the claim follows by Lemma \ref{GoF}.
\end{Pf}
We also need a variation of Corollary \ref{Kanten}.
\begin{Not}\label{BNot}
In the situation of Notation \ref{KantenNot}, let 
$$B:= {(C\to d)\times\II}\push{(C\to d)\times\{1\}}{(C\to d')}.$$
Let 
$$\map{r_B}{B}{(C\to d')}$$
be the projection onto the first factor on $(C\to d)\times\II$ and the 
identity on the rest. Let $l_B$ be the left adjoint to $r_B$, i.e. $l_B$
sends $(C\to d)$ to $(C\to d)\times\{0\}$. Let further
$$\map{p_B}{B}{\II}$$
send $(C\to d')$ to $1$ and $(C\to d)\times\{0\}$ to $0$, and let
$$j_B:=j_{d'}\circ r_B.$$
\end{Not}
\begin{Cor}\label{Variation}
With Notation \ref{BNot}, we have a functor isomorphism
$$(d\leq d')^*\circ\holkan{f}\cong\holkan{p_B}j_B^*.$$
\end{Cor}
\begin{Pf}{}
We have
$$p^{d'}_d = p_B\circ l_B.$$
Therefore
$$\holkan{p^{d'}_d}j^*_{d'}\cong\holkan{p_B}\holkan{l_B}j^*_{d'}
\cong\holkan{p_B}r_B^*j^*_{d'}.$$
\end{Pf}
Next we discuss the compatibility of the smash product with the
triangulated structure. 
\begin{Prop}\label{smcone}
There is a functor isomorphism
$$\cone(-)\wedge^L\cone(-)\cong\cone(-\Box^L-).$$
\end{Prop}
\begin{Pf}{}
We have 
\begin{eqnarray*}
\cone (f)\wedge^L\cone (g) & = & \hocolim{\Vv}\circ\Vst (f)\wedge^L
                                 \hocolim{\Vv}\circ\Vst (g)\\
&\cong& \hocolim{\Vv\times\Vv}
                                   (\Vst (f) \wedge^L\Vst (g)),\\
\end{eqnarray*}
where the first isomorphism is Definition \ref{conedef} and the 
second isomorphism is Corollary \ref{smhol}.
The diagram
$$
\left(\Vst f\right) \wedge^L\left(\Vst g\right)
$$ 
has the form (for $f\negmedspace : X\to Y$ and $g\negmedspace : U\to V$)
$$
{\Setunitlength
\begin{picture}(12,8)
\def\Kr{\circle{.4}}
\put(0.8,6.8){$\star$}
\put(4.8,6.8){$\star$}
\put(2.8,5.8){$\star$}
\put(3.8,2.8){$\star$}
\put(6,2){\Kr}
\put(6.8,6.8){$\star$}
\put(8,3){\Kr}
\put(9,6){\Kr}
\put(11,7){\Kr}  
\put(6.1789,2.0894){\line(2,1){1.6422}}
\put(5.8211,2.0894){\line(-2,1){1.6422}}
\put(3.1789,6.0894){\line(2,1){1.6422}}
\put(2.8211,6.0894){\line(-2,1){1.6422}}
\put(9.1789,6.0894){\line(2,1){1.6422}}
\put(8.8211,6.0894){\line(-2,1){1.6422}}
\put(4.12,3.16){\line(3,4){2.76}}
\put(3.88,3.16){\line(-3,4){2.76}}
\put(6.12,2.16){\line(3,4){2.76}}
\put(5.88,2.16){\line(-3,4){2.76}}
\put(8.12,3.16){\line(3,4){2.76}}
\put(7.88,3.16){\line(-3,4){2.76}}
\put(6,1){\makebox(0,0)[ca]{\small{$X\wedge U$}}}
\put(8.5,2.5){\makebox(0,0)[la]{\small{$Y\wedge U$}}}
\put(9.5,5.5){\makebox(0,0)[la]{\small{$X\wedge V$}}}
\put(13,7){\makebox(0,0)[cb]{\small{$Y\wedge V$}}}
\end{picture}
.}
$$
In particular, Corollary \ref{smzur} implies
$$\left(\Vst(f)\wedge^L\Vst(g)\right)
\at{((0,0)\leq(1,0))\times((0,0)\leq(0,1))}
\cong f\wedge^L g \in\Ho(\MM^{\Ii\times\Ii}).$$
We consider the map
\begin{eqnarray*}
\specialmap{pr} \VVv\times\VVv & \to & \VVv \\
((0,1),(0,1)) & \mapsto  & (0,1) \\
((0,1),(0,0)),((0,0),(0,0)),((0,0),(0,1)) & \mapsto & (0,0)\\
(1,0)\times\VVv\cup\VVv\times (1,0) & \mapsto & (1,0).
\end{eqnarray*}
By Corollary \ref{Kanten} and with $\pC$ as in (\ref{box}), we have 
\begin{eqnarray*}
\holkan{pr}\left(\Vst f\wedge^L\Vst g\right)\at{(0,0)\leq(0,1)} & \cong &
                         \holkan{\pC}\left(f\wedge^L g\right) \\
&  \cong & f \Box^L g ,
\end{eqnarray*}
where the second isomorphism is (\ref{lbox}).
Further, by (\ref{ecken}), we have
$$\holkan{pr}\left(\Vst f \wedge^L\Vst g\right)\at{(1,0)} \cong \star .$$
Together, we obtain
$$\holkan{pr}\left(\Vst f\wedge^L\Vst g\right) \cong \Vst\left(f\Box^Lg\right) 
$$
If we apply $\hocolim{\Vv}$ to this equation, the claim follows.
\end{Pf}
For later reference, we recall another functor isomorphism 
\cite[Thm.2]{Franke:95}:
We have
\begin{equation}\label{holCone}
\Cone\circ\holkan{pr_\Ii}(-) \cong \holkan{pr_\Ii}\circ\Cone_C(-),
\end{equation}
where 
$$\map{pr_\Ii}{C\times\II}{\II}$$
denotes the projection to the second factor.
It follows that
\begin{equation}\label{holcone}
\cone\circ\holkan{pr_\Ii} (-)\cong \hocolim{C}\circ\cone_C(-).
\end{equation}
The remainder of this section is about the interaction of $\Cone$ with
the monoidal structure. It needs some preparation and is not needed in
the proof of Theorem \ref{Maintheorem}. The reader only interested in
this theorem can skip ahead to Section \ref{welchesmodell}.
Before we can discuss the compatibility of the smash product with the
"cone inclusion" functor $\Cone$ from (\ref{Conedef}), we
need an alternative description of $\Cone$. For completeness, we
also give a similar description of the "cone map"
$$\map{\Cone(\Cone(f))}{\cone(f)}{\Sigma X}\in\Ho(\MI).$$
For our next definition, we make an exception from our conventions
about posets, and let the arrows point to the right and down.
\begin{Def}
We define the functors
$$\map{\kegink, \kegabb}{\MI}{\MM^{\Ii\times\Ii}},$$
that map
$$\xymatrix{X\ar[r]^f&Y}\in\MI$$
to
$$
\xymatrix{\star\ar[r]\ar[d] & Y\ar@{=}[d]\\
X\ar[r]^f&Y}\raisebox{-2em}{\text{  and  }}
\xymatrix{X\ar[r]^f\ar@{=}[d] & Y\ar[d]\\
X\ar[r]& \star}
\raisebox{-2em}{$\in\MM^{\Ii\times\Ii}$}
$$
respectively.
\end{Def}
In these pictures, horizontal arrows correspond to the second factor of
$\II\times\II$ whereas vertical arrows correspond to the first factor.
Both of these functors preserve weak equivalences and therefore induce
functors on the homotopy categories. We use the same names for
these induced functors.
\begin{Lem}\label{Kegelabb}
There are functor isomorphisms
$$\Cone(-)\cong\cone_\Ii\circ\kegink(-)$$
and
$$
\Cone\circ\Cone(-)\cong\cone_\Ii\circ\kegabb(-).
$$
\end{Lem}
Here we used the following notation: Let $C$ be a poset. We write
$$\map{\cone_C}{\Ho(\MM^{C\times\Ii})=\Ho((\MM^C)^\Ii)}{\Ho(\MM^C)}
$$
for the cone functor with $\MM^C$ in the role of $\MM$. 
In our situation, $C = \II$, the first factor of $\II\times\II$.

\begin{Pf}{}
We write
${\Vst_C}$ for the functor $\Vst$ with $\MM^C$ in the role of $\MM$.
We claim that
$$\holkan{\id_\Ii\times(\Vv\subset(\Ii\times\Ii))}\circ\Vst_\Ii\circ\kegink(
X\stackrel{f}{\to}Y)$$
is of the form
$$
\xymatrix{ & \star\ar[rd] &\\
\star\ar@{=}[ru]\ar[r]\ar[d] & Y \ar@{=}[d]\ar@{=}[r] & 
Y\ar[d]^{\cone_\Ii(\kegink(f))}\\
X\ar[r]^f\ar[rd] & Y\ar[r]^{\Cone(f)} & \cone(f).\\
& \star\ar[ru]}
$$
But this follows from Definition \ref{conedef}, (\ref{Conedef}), and
the fact that for a diagram
$X\in\Ho(\MM^{C\times\Ii})$ 
the vertex
$\left(\cone_C(X)\right)_c$
is isomorphic to the cone of the corresponding restriction $X\at{c\times\Ii}$ (see \cite[1.4.2]{Franke:95}).

For the cone map, we look at
$$\holkan{\id_\Ii\times(\Vv\subset(\Ii\times\Ii))}\circ\Vst_\Ii\circ\kegabb(
\xymatrix{X\ar[r]^f&Y}).$$
By the same argument as above and by (\ref{ecken}) this is of the shape
$$
\xymatrix{& \star\ar[rd] \\
X\ar[r]_f\ar[ru]\ar@{=}[d] & Y \ar[r]_{\Cone(f)}\ar[d] & \cone(f) 
\ar[d]^{\cone_\Ii(\kegabb(f))} \\
X\ar[r]\ar[rd] & \star \ar[r] & \Sigma X \\
& \star.\ar[ru]
}
$$
The right vertical edge is given by (\ref{ecken}) with $\MI$ in the
role of $\MM$. It is
$$\cone_\Ii\circ\kegabb(f).$$
We have to show that the right square is homotopy bi-cartesian.
But the top square is homotopy bi-cartesian, and so is the square
that we obtain by putting the top square and the right square next to
each other.
Therefore, \cite[Prop.1.4.6]{Franke:95} implies that the right square
is also homotopy bi-cartesian.
\end{Pf}
Let now $\mathcal N$ stand for the model category $\MI$ with vertex-wise
cofibrations and weak equivalences.
Then
$$\map{\kegink}{\MI}{\mathcal N^\Ii}$$ 
preserves Reedy cofibrations.
Here we identified the first exponent in $\mathcal N^\Ii = (\MI)^\Ii$
with the vertical arrows.
Therefore Lemma \ref{GoF} implies
\begin{Cor}
There is a functor isomorphism
$$\kegink\circ(-\Box^L -)\cong\left(\kegink(-)\right)\Box^L_\Ii
\left(\kegink(-)\right).$$
\end{Cor}
Here $\Box_\Ii$ denotes the pushout with $\mathcal N$ in the role of
$\MM$. 
Proposition \ref{smcone} now implies
\begin{Cor}\label{smCone}
There is a functor isomorphism
$$\Cone(-\Box^L-)\cong\Cone(-)\wedge^L_\Ii\Cone(-).$$
\end{Cor}
Here $\wedge_\Ii$ denotes the (internal) smash product in $\mathcal N$.

\begin{Pf}{}
We have
\begin{eqnarray*}
\Cone(f\Box^Lg) & \cong & \cone_\Ii\left(\kegink(f\Box^Lg)\right)\\
&\cong & \cone_\Ii\left(\left(\kegink(f)\right)
\Box^L_\Ii\left(\kegink(g)\right)\right)\\
&\cong&\cone_\Ii\left(\kegink(f)\right)
\wedge^L_\Ii\cone_\Ii\left(\kegink(g)\right)\\    
&\cong&\Cone(f)\wedge^L_\Ii\Cone(g),
\end{eqnarray*}
where the third isomorphism is Proposition \ref{smcone}.
\end{Pf}
\section{Which model?}\label{welchesmodell}
There are four canonical choices of model for $\SE$, all of which
are equally well suited for our purposes. Firstly, we can work either in
the world of
symmetric spectra \cite{HSS} or in the world of $S$-modules \cite{EKMM}. 
Let $\MM_S$ be one of these two models for the stable homotopy category. 
One possible model for $\SE$ is $\MM_S\lbrack E(1)^{-1}\rbrack$ (compare
Example \ref{BL}). To obtain the other model, we recall that Hopkins and
Ravenel have shown that localization at $E(n)$ is smashing, i.e. that
for all $X$ in $\mathcal S$, one has
$$X_{E(n)}\cong X\wedge_S^L S_{E(n)}$$
\cite{Ravenel:92}. 
If the localization at a spectrum $E$ is smashing, the fibrant
replacement $S_E$ of the sphere spectrum in $\MM_S\lbrack E^{-1}\rbrack$ can
be chosen to be a strict, strictly commutative ring spectrum, and 
$$\map{-\wedge_SS_E}{\MM_S\lbrack E^{-1}\rbrack}{\MM_{S_E}}$$
is a Quillen equivalence. In other words,
$$\Ho(\MM_{S_E})$$ 
is the localization of $\MM_S$ at $E_*$, and 
$$-\wedge_S^LS_E$$
is the localization 
functor\footnote{For symmetric spectra, this statement is
  \cite[3.2.(iii)]{SS2}. 
  For $S$-modules this is a result of Wolbert
  \cite{Wolbert}, which can also be found in \cite{EKMM}. In order to
  get this precise statement from \cite{EKMM}, one actually has to combine 
  a few propositions: it follows from VIII.2.1.,
  VIII.3.2., and from the fact that by III.4.2. and VII.4.9. the
  derived categories of $\Lambda S$-modules and of $S$-modules
  are equivalent as monoidal categories. Here $\Lambda S\to S$ denotes
  the q-cofibrant replacement of the sphere spectrum.}.
For Franke's methods it is irrelevant which model for $\SE$ one likes
to choose. The only property of the model category that is relevant
for him is that it induces a system of {\it triangulated} diagram
categories, i.e. that the model category is 
stable\footnote{A model category is called stable, if the suspension functor
  is invertible in the homotopy category. In our example this follows
  from the fact that the suspension in $\Ho(\MM_R)$ is given by
  smashing over $R$ with $S^1\wedge_S^L R$. Therefore smashing over
  $R$ with $S^{-1}\wedge_S^L R$ is an inverse of the suspension.
  Also Quillen equivalent models give rise to the same suspension
  functor. 
  The fact that the homotopy category of a stable model category is
  triangulated, is proved in \cite[7.1]{Hovey}.}. 
Moreover, Schwede \cite{Schwede} has constructed a 
functor\footnote{just for the moment, $\MM_S$ denotes the $S$-modules from
  \cite{EKMM} and $Sp^\Sigma$ denotes the symmetric spectra from
  \cite{HSS}.}
$$\map{\Phi}{\MM_S}{Sp^\Sigma},$$
that maps (strict) ring spectra to (strict) ring spectra, and
induces for any strict (strictly commutative) ring spectrum $R$ a
monoidal equivalence
$$\Ho(\MM_R)\stackrel{\sim}{\longrightarrow}\Ho(\operatorname{\Phi(R)-mod}).$$

We also know that for any strict, strictly commutative ring spectrum
$R$, the functor
$$\map{-\wedge_SR}{(\MM_S,\wedge_S)}{(\MM_R,\wedge_R)}$$
is strictly monoidal. It follows that up to equivalence 
all four models mentioned
above give rise to the same system of triangulated diagram
categories, and the same smash product on it.
\section{The Derived Tensor Product}
\label{flatmodel}
In this section we define the derived tensor product
$\otimes^L_{E(1)_*}$ on the derived category of quasi-periodic 
cochain complexes in $\eeO$. 
We show that $\otimes^L_{\E}$ is a monoidal structure. 
Since $\operatorname{Comod_{\EE}}$ does not have enough projectives, we
need to work with flat replacement rather than projective replacement.
Our definition of a flat complex is object-wise, forcing us to replace
both sides of $\otimes_{\E}$ with flat objects.
%
%
%
%
The following lemma is essentially due to Christensen and Hovey.
\begin{Lem}\label{Hovey-Lem}
Let $E$ be a Landweber exact cohomology theory, and let $C$ be a
quasi-periodic cochain complex of $(E_*, E_*E)$-comodules. Then there exist a
flat quasi-periodic cochain complex $P$ of $(E_*, E_*E)$-comodules
and a quasi-periodic quasi-isomorphism $i$ from $P$ to $C$.
Moreover, they pair $(P,i)$ depends on $C$ in a functorial way.
\end{Lem}
\begin{Pf}{}
  We claim that   
  cofibrant replacement in Christensen and Hovey's projective model
  structure\footnote{Cf.\ \cite[2]{Hovey:projective} and 
    \cite[4.4]{Christensen:Hovey}.}
  on $\mathcal Ch(E_*,E_*E)$  can be done in such a way that these
  conditions are satisfied.
  We have to show three things: (a) Cofibrant objects in the
  projective model structure are flat over $E_*$, (b) weak equivalences in the
  projective model structure are quasi-isomorphisms, and (c) if $C$ is
  quasi-periodic, one can
  choose a quasi-periodic cofibrant replacement of $C$. 
  \begin{enumerate}
    \renewcommand{\labelenumi}{(\alph{enumi})}
    \item
      Let $P$ be cofibrant in the projective model structure. 
      By \cite[4.4]{Christensen:Hovey},
      $P$ is a retract of a (transfinite) colimit of a
      diagram of complexes
      $$
        0=P_0\to P_1 \to\cdots\to P_\alpha\to P_{\alpha+1}\to\cdots,
      $$ 
%
%
      where each $P_\alpha\to P_{\alpha+1}$ is a degree-wise split
      monomorphism whose cokernel is a complex of so called ``relative
      projectives''  
      with no 
      differential, and if $\alpha$ is a limit ordinal, $P_\alpha$ is
      the colimit over all $P_\beta$ with $\beta<\alpha$.
      It is also pointed out in \cite[p.15]{Hovey:projective} that every
      ``relative projective'' comodule is 
      projective as an $E_*$-module.
      We can therefore prove by transfinite induction that the limit
      over all the $P_\alpha$ is flat: The complex $P_0$ is flat.
      In every degree, $P_{\alpha+1}$ is the direct sum of
      $P_\alpha$ with a projective $E_*$-module, thus if $P_\alpha$ is
      degree-wise flat, so is $P_{\alpha +1}$. If $\alpha$ is a limit
      ordinal, $P_\alpha$ is a direct limit of flat objects and hence flat.
      For the same reason the colimit over the entire diagram
      is flat. 
      Moreover, in an abelian category, retracts are direct summands, and
      thus retracts of flat modules are flat. This proves that $P$ is
      degree-wise flat.
    \item
      Since $E_*$ is Landweber exact, it follows from
      \cite[Sec.1.4]{Hovey:projective} that $(E_*,E_*E)$ satisfies the
      conditions of \cite[2.1.5]{Hovey:projective}. 
      Therefore, weak equivalences are quasi-isomorphisms.
    \item
      Let now $C$ be quasi-periodic. In the proof of
      \cite[4.2]{Christensen:Hovey} an explicit cofibrant replacement
      is constructed. Observe that for a quasi-periodic complex $C$ the
      complexes $P_i$ and $Q_i$ may be chosen in such a way that the
      ``partial cofibrant replacement'' of $C$ is again quasi-periodic.
      None of the other steps in \cite[4.2]{Christensen:Hovey}
      (colimits, 
      pullbacks, path objects, cofibres) destroy quasi-periodicity.
      This proves the claim.  
  \end{enumerate}
\end{Pf}
\begin{Rem}\label{flat-Rem}
Flatness in our case means flatness over $\E$.
Therefore flat objects are exactly the $p$-torsion free
objects, and flat replacement in $\mathcal C^1(\eee)$
translates into flat replacement in $\mathcal C^{2p-2}(\eeO)$.
Note also that for the same reason subobjects of flat objects are
again flat.
\end{Rem}
The following is a corollary of \cite[III.2.10]{Gelfand:Manin}. 
\begin{Prop}
  Let $\mathcal C$ be a category, $S$ a left-localizing system of
  morphisms in $\mathcal C$ and $\mathcal B\subseteq\mathcal C$ a full
  subcategory. Suppose that for all $X\in\Ob\mathcal C$ there exist 
  $Y\in\Ob\mathcal B$ and $s\negmedspace :Y\to X$ in $S$. Then 
  $$
    S_{\mathcal B} := S\cap\Mor\mathcal B
  $$
  is left-localizing in $\mathcal B$, and the canonical map
  $$
    {\mathcal B}[S_{\mathcal B}^{-1}]\to\mathcal C[S^{-1}]
  $$
  is an equivalence of categories.
\end{Prop}
\begin{Not}
  Let $$\mathcal C_{flat}\subset \mathcal C^{2p-2}(\eeO)$$
  be the full subcategory of flat objects, and let
  $$\map{\gamma_{flat}}{\mathcal C_{flat}}{\mathcal D^{2p-2}(\eeO)}$$
  be the composition of its inclusion with the localization functor.
  Let $K_{flat}$ denote the homotopic category 
  (in the sense of homological algebra) of $\mathcal C_{flat}$, and
  let $S$ denote the class of quasi-isomorphisms in  $K_{flat}$.
\end{Not}
\begin{Cor}\label{flat-equiv-Cor}
  The functor $\gamma_{flat}$ induces an equivalence of localized
  categories  
  $$
    K_{flat} [S^{-1}] \to \mathcal D^{2p-2}(\eeO).
  $$
\end{Cor}
In other words, flat replacement can be done in such a way that it is
functorial on morphisms in the derived category, not just the strict
category: let $\phi$ be an
equivalence of categories inverse to the equivalence in
Corollary \ref{flat-equiv-Cor}, then $\phi$ is such a flat
replacement functor on the derived category.
\begin{Cor}
  The tensor product in $\mathcal C^{2p-2}(\eeO)$ has a left derived
  functor $\otimes^L_{\E}$, defining a symmetric monoidal structure on the
  derived category.
\end{Cor}
\begin{Pf}{}
  The functor
  $$
    \map{(-\otimes_{\E}-)}{K_{flat}\times
    K_{flat}}{K_{flat}[S^{-1}]} 
  $$
  takes pairs of acyclic complexes to acyclic complexes. Therefore, by
  the same argument as in \cite[III.2.23]{Gelfand:Manin},
  it factors over the localization 
  $$
    (K_{flat}\times K_{flat})[(S\times S)^{-1}] =
    K_{flat}[S^{-1}]\times K_{flat}[S^{-1}].
  $$
  As in the proof of \cite[III.6.8]{Gelfand:Manin}, we precompose
  $(-\otimes_{\E}-)$ with $\phi\times\phi$, 
  to obtain a functor
  $(-\otimes^L_{\E}-)$ satisfying the universal property of the left
  derived functor of $(-\otimes_{\E}-)$. As in the proof of
  \cite[4.3.2]{Hovey}) 
  the structure diagrams making $\otimes_{\E}$ a monoidal structure 
  can be translated into the respective diagrams for $\otimes_{\E}^L$.
\end{Pf}
\section{Proof of the main theorem}\label{beweisvonthm}
In this section we prove Theorem \ref{Maintheorem}. The
result of Section \ref{flatmodel} allows us to restrict our discussion
to flat complexes.
We need to show the following.
\begin{Thm}\label{grossersatzvonganter}
There is a functor isomorphism
$$(-\wedge^L_{S_{E(1)}}-)\circ\Rek\circ\gamma_{flat}\cong
\Rek\circ\gamma_{flat}\circ(-\otimes_{E(1)_*}-).$$
\end{Thm}
Let $C$ and $\tilde C$ be in $\mathcal C_{flat}$.
Recall from (\ref{Recdef}) that
$$\Rek(C) = \hocolim{C_N} Q^{-1}(C),$$
where $Q^{-1}(C)$ is the object of $\mathcal L$
defined in Notation \ref{QNot}, and $\mathcal L$ is as in Notation
\ref{LNot}. We need to identify
$$\Rek^{-1}(\Rek(C)\wedge^L_{S_{E(1)}}\Rek(\tilde C))$$
with $C\otimes_{\E}^L\tilde C$.
By (\ref{Recdef}), there are two steps to computing the left
hand side: First, we need to find an object of $\mathcal L$ whose
homotopy colimit is 
$$\Rek(C)\wedge^L_{S_{E(1)}}\Rek(\tilde C).$$
We do this in \ref{firststep}. Then, in \ref{secondstep}, we apply $Q$
to this object.
In the following, most computations take place on the homotopy 
level, and we allow ourselves the following abuse of notation:
\begin{Not}
We abbreviate
$-\wedge^L_{S_{E(1)}}-$ by $\quad -\wedge -$.
\end{Not}
\subsection{An object of $\mathcal L$, whose homotopy
colimit is $\Rek(C)\wedge\Rek(\tilde C)$}\label{firststep}
\begin{Not}\label{ANot}
We write $A$ for $Q^{-1}(C)$ and $\tilde A$ for $Q^{-1}(\tilde C)$.
\end{Not}
We have
$$\Rek(C)\wedge\Rek(\tilde C)\cong(\hocolim{C_N}A)\wedge(\hocolim{C_N}\tilde A)
\cong \hocolim{C_N\times C_N}(A\wedge\tilde A),$$
where the first isomorphism is (\ref{Recdef})
and the second is Corollary \ref{smhol}.
Thus
$A\wedge\tilde A$ is a diagram with the correct homotopy colimit. Of
course it is not an object of $\mathcal L$, it does not even have the
right shape.
A way to change the shape of a diagram without changing the homotopy
colimit is to apply a homotopy left Kan extension along a map of
posets to it. 
We need to pick this map out of $C_N\times C_N$
such that the $\E$-homology of the vertices remains
concentrated in the correct degrees. Since there is no such map to
$C_N$, we choose as target the poset
$$D_N := \{\beta_n, \gamma_n, \zeta_n\mid n\in\zn\}$$ 
with relations generated by
$$\beta_{n+1}\leq\gamma_n,\quad \beta_n\leq\gamma_n,
\quad \gamma_{n+1}\leq\beta_n
\text{ and } \gamma_n\leq\zeta_n.
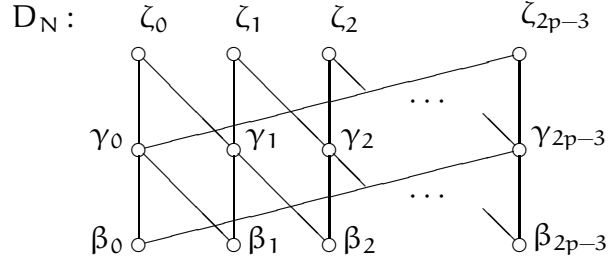
$$ 
\begin{center}
$$\posetDN$$
\end{center}
The map of posets is
\begin{eqnarray*}
pr\negmedspace : C_N\times C_N & \to    & {D_N}\\
(\beta_s, \beta_t)            &\mapsto & \beta_{s+t} \\ 
(\beta_s, \zeta_t),(\zeta_s, \beta_t) &\mapsto &\gamma_{s+t} \\
(\zeta_s, \zeta_t) & \mapsto & \zeta_{s+t}.
\end{eqnarray*}
\begin{Not}\label{ENot}
With $A$ and $\tilde A$ as in Notation \ref{ANot}, let
$$
E := \holkan{pr}(A\wedge\tilde A).
$$
\end{Not}
This is an object of $\Ho(\EM^{D_N})$.
\begin{Prop}\label{inj}
The objects
$$\E(E_{\alpha_n})\text{ with }\alpha\in\{\beta,\gamma,\zeta\}$$ 
are concentrated in degrees congruent to $n$ modulo $2p-2$.
The morphisms
\begin{equation}
E(1)_*(E_{\gamma_n})\to E(1)_*(E_{\zeta_n}),
\end{equation} 
induced by the corresponding edges in $E$, are monomorphisms.
\end{Prop}
\begin{Pf}{}
By Corollary \ref{Kanten}, 
\begin{equation}\label{gammantozetan}
E\arrowvert_{\gamma_n\leq\zeta_n} \cong 
\holkan{p^{\zeta_n}_{\gamma_n}}j_{\zeta_n}^*(A\wedge\tilde A)
\end{equation}
(Notation \ref{KantenNot}). Still in Notation \ref{KantenNot},
%
%
%

\noindent $\CA{\zeta_n} = $
\begin{center}
  $$ 
    \nbutterfly{(\zeta_s,\zeta_t)}{(\zeta_s,\beta_t)}
    {(\beta_s,\zeta_t)}{(\beta_s,\beta_t)}{(\beta_{s+1},\zeta_t)}
    {(\zeta_s,\beta_{t+1})}{(\beta_{s+1},\beta_t)}{(\beta_s,\beta_{t+1})}
    {(\beta_{s+1},\beta_{t+1})}
  $$
\end{center}
The vertices marked black are the ones mapped
to zero by $p_{\gamma_n}^{\zeta_n}$. In other words, 
the source of (\ref{gammantozetan}) is
the homotopy
colimit over the black subdiagram, whereas the target is the homotopy 
colimit over the entire diagram.
We need to compute the $\E$-homology of (\ref{gammantozetan}).
We do so by applying spectral sequence (\ref{SF}) to the right
hand side of (\ref{gammantozetan}).
The $E_2$-term involves left derived Kan extensions along 
$p_{\gamma_n}^{\zeta_n}$, and the following three lemmas are about
their computation.
\begin{Not}\label{VO}
Let 
$$\VO_n\subset (\CCto\zeta_n)$$
denote the sub-poset with elements
$$\{(\zeta_s,\zeta_t),(\zeta_s,\beta_t),(\beta_s,\zeta_t),(\beta_s,\beta_t),
(\beta_{s+1},\beta_t)\mid s+t=n\}.$$ 
Let
$$\map\jo{\VO_n}{(\CCto\zeta_n)},$$
denote its inclusion, and
$$\map{\lo}{(\CCto\zeta_n)}{\VO_n}$$
denote the left adjoint to $\jo$.
Let further
$$\map{\po:=p_{\gamma_n}^{\zeta_n}\circ\jo}{\VO_n}{\II}.$$
\end{Not}
\begin{Lem}
With this notation we have functor isomorphisms
$$(\lkan{p_{\gamma_n}^{\zeta_n}})_s =(\lkan{\po})_s\circ\jo^*.$$
\end{Lem}
\begin{Pf}{}
By adjointness,
$$\lkan{\lo}=\jo^*,$$
and the higher derived left Kan extensions along $\lo$ vanish. We have
$$p_{\gamma_n}^{\zeta_n}=\po\circ\lo,$$
and the Grothendieck spectral sequence for the composition of derived
functors proves the claim.
\end{Pf}
\begin{Not}
Consider the poset
$$\VY_n := \{\alpha_{s,t},(\zeta_s,\zeta_t),(\beta_{s+1},\beta_t)\mid s+t=n\}$$
with relations
\begin{eqnarray*}
(\beta_{s+1},\beta_t),(\beta_s,\beta_{t+1}) & \leq & \alpha_{s,t}\\
\alpha_{s,t}                                & \leq & (\zeta_s,\zeta_t).
\end{eqnarray*}
Let $g$ and $\py$ be the maps of posets
\begin{eqnarray*}
g\negmedspace : \VO_n                      & \to     & \VY_n              \\
(\zeta_s,\zeta_t)                          & \mapsto & (\zeta_s,\zeta_t)  \\
(\zeta_s,\beta_t),(\beta_s,\zeta_t),(\beta_s,\beta_t)& \mapsto & \alpha_{s,t}\\
(\beta_{s+1},\beta_t)                      & \mapsto &(\beta_{s+1},\beta_t)\\
\end{eqnarray*}
and
\begin{eqnarray*}
\py\negmedspace : \VY_n                      & \to     & \II              \\
(\zeta_s,\zeta_t)                          & \mapsto & 1  \\
(\beta_{s+1},\beta_t),\alpha_{s,t}         & \mapsto & 0.  \\
\end{eqnarray*}
\end{Not}
\begin{Lem}\label{hilfsatz2}
Let $X\in\Ho(\MM^\Vo)$ be such that
$$X\at{\{(\zeta_s,\beta_t),(\beta_s,\zeta_t),(\beta_s,\beta_t)\}}$$
is $\colim{\Vv}$-acyclic. Then there is an isomorphism (natural in
such $X$)
$$(\lkan\po)_sX = (\lkan\py)_s\circ\lkan gX.$$
\end{Lem}
\begin{Pf}{}
We have 
$$p_\Vo = p_\Vy\circ g.$$
By (\ref{striktekanten}), $X$ is $\lkan{g}$-acyclic. The Grothendieck
spectral sequence implies the claim.
\end{Pf}
Consider a diagram $\mathcal Y$ with underlying poset $\VY$:
\begin{equation}\label{Ypsilon}
\xymatrix@=2ex{
& Z_1 && Z_2 && Z_3 &\\
& X_1\ar[u]_{h_1} && X_2\ar[u]_{h_2} && X_3\ar[u]_{h_3} &\\
Y_1\ar[ru]^{\dots \quad g_1} && Y_2\ar[lu]^{f_1}\ar[ru]^{g_2}
&&Y_3\ar[lu]^{f_2}\ar[ru]^{g_3} && Y_4\ar[lu]_{f_3\quad\cdots}\\
}
\end{equation}
\begin{Lem}\label{hilfsatz3}
Let $\mathcal Y$ be as in (\ref{Ypsilon}).
Then there are functor isomorphisms
$$\lkan{p_\Vy}(\mathcal Y)\cong
\coeq_\Ii \left({\Setunitlength\raisebox{2\unitlength}{
\xymatrix{\bigoplus Y_i \ar @<.6ex> [r] ^{\oplus h_if_i}
   \ar @<-.6ex>  [r] _{\oplus h_ig_i} & \bigoplus Z_i\\
\bigoplus Y_i \ar @<.6ex> [r] ^{\oplus f_i}
   \ar @<-.6ex>  [r] _{\oplus g_i} \ar@{=}[u]
& \bigoplus X_i\ar[u]_{\oplus h_i}\\
}}}
\right),
$$
$$(\lkan{p_\Vy})_1(\mathcal Y)\cong
\equalizer_\Ii \left({\Setunitlength\raisebox{2\unitlength}{
\xymatrix{\bigoplus Y_i \ar @<.6ex> [r] ^{\oplus h_if_i}
   \ar @<-.6ex>  [r] _{\oplus h_ig_i} & \bigoplus Z_i\\
\bigoplus Y_i \ar @<.6ex> [r] ^{\oplus f_i}
   \ar @<-.6ex>  [r] _{\oplus g_i} \ar@{=}[u]
& \bigoplus X_i\ar[u]_{\oplus h_i}\\
}}}
\right),
$$
where $\operatorname{(co)eq}_\Ii$ denotes the (co)equalizers of the
pairs of horizontal arrows (together with the induced map between
them). The higher derived left Kan extensions along $p_\Vy$ vanish.
\end{Lem}
\begin{Pf}{}
The first statement is (\ref{striktekanten}), the third statement
follows from (\ref{striktekanten}) and the fact that for a diagram of
length one the second and all higher derived colimits vanish. The
second statement follows from the other two together with the
universal property of derived functors and the snake lemma applied to
an $\lkan{}$-acyclic resolution.
\end{Pf}
We are now ready to complete the proof of Proposition \ref{inj}.
As subobjects of flat objects,
$\E(A_{\alpha_s})$ and $E(1)_*(\tilde A_{\alpha_t'})$ are flat, and 
the Künneth spectral 
sequence\footnote{The Künneth spectral sequence for non connective 
spectra can for example be found in \cite{EKMM}.} 
reduces to
$$\E (A_{\alpha_s}\wedge\tilde A_{\alpha_t'}) 
\cong \E (A_{\alpha_s})\otimes_{\E}\E (\tilde
A_{\alpha_t'}).$$ 
Firstly, we have

\noindent
$E(1)_{-n}(j_{\zeta_n}^*\AAA)=$
\nopagebreak
\begin{center}
\begin{equation}\label{erstesbild}
\nbutterfly{Z^s\otimes\tilde Z^t}
{Z^s\otimes\tilde B^t}{B^s\otimes\tilde Z^t}{B^s\otimes\tilde B^t}{0}{0}{0}{0.}
{0}
\end{equation}
\end{center}
Since $B^s$ and $\tilde B^t$ are flat, the maps
$B^s\otimes\tilde B^t\to B^s\otimes\tilde Z^t$ and
$B^s\otimes\tilde B^t\to Z^s\otimes\tilde B^t$ are monomorphisms.
Thus $j_\Vo$ applied to (\ref{erstesbild}) satisfies the condition of
Lemma \ref{hilfsatz2}. We are reduced to the situation of Lemma
\ref{hilfsatz3}, with 
$$\mathcal Y = \lkan{g}j_\Vo^*((\ref{erstesbild})).$$
One computes $\mathcal Y$ using (\ref{striktekanten}).
Together with Lemma \ref{hilfsatz3}, we obtain
%
$$\lkan{p_{\gamma_n}^{\zeta_n}}E(1)_{-n}(j_{\zeta_n}^*\AAA) = 
\left (\Oplus_{s+t=n}Z^s\otimes\tilde B^t
\push{B^s\otimes\tilde B^t}B^s\otimes\tilde Z^t\hookrightarrow
\Oplus_{s+t=n}Z^s\otimes\tilde Z^t\right ),$$
the canonical inclusion, and
the higher derived left Kan extensions of (\ref{erstesbild}) 
along $p_{\gamma_n}^{\zeta_n}$ vanish.

Secondly, we consider

\pagebreak
\noindent
$E(1)_{-n-1}(j^*_{\zeta_n}\AAA)=$
\nopagebreak
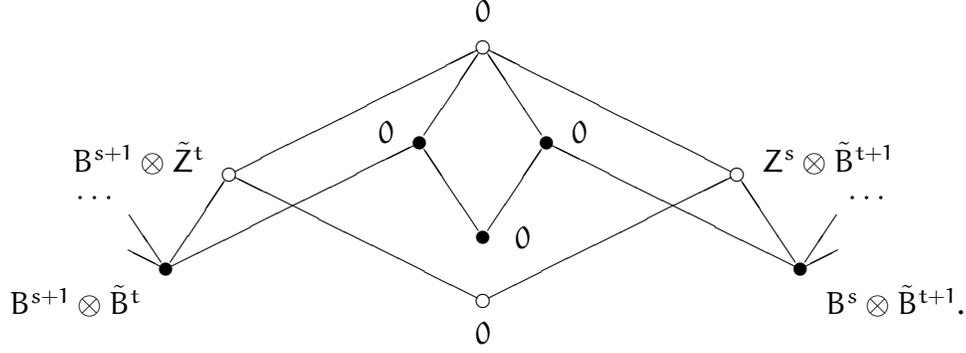
\begin{center}
$$
\nbutterfly{0}{0}{0}{0}{B^{s+1}\otimes\tilde Z^t}
{Z^{s}\otimes\tilde B^{t+1}}{B^{s+1}\otimes\tilde B^t}
{B^{s}\otimes\tilde B^{t+1}.}{0}
$$\end{center}
The argument is as above, the result is:
\begin{eqnarray*}
\lkan{p_{\gamma_n}^{\zeta_n}}E(1)_{-n-1}(j^*_{\zeta_n}\AAA) 
           &=& 0\negmedspace : 0\to 0,\\
\lkan{p_{\gamma_n}^{\zeta_n}}{}_1E(1)_{-n-1}(j^*_{\zeta_n}\AAA) &=&
\id_{\Oplus_{s+t=n+1}B^s\otimes\tilde B^t}.
\end{eqnarray*}
The higher derived left Kan extensions vanish.

Thirdly, 
$$E(1)_{-n-k}(\AAA\at\Vo) = 0,$$ 
if $k$ is not congruent to $0$ or $1$ modulo $2p-2$.
This completes the calculation of the $E_2$-term.
It is concentrated
in degrees $(0,m)$ and $(-1,m+1)$, with $m\equiv n\hbox{$\mod 2p-2$}$.
Therefore the spectral sequence collapses at the $E_2$-term and
becomes a short exact sequence (of morphisms)
\begin{equation}\label{BZ}
\xymatrix@=3ex{
0\ar[r]&{\Oplus_{s+t=n}Z^s\otimes\tilde Z^t\ \ }
\ar[0,1]&
{E(1)_{-n}(E_{\zeta_n})}
\ar[0,1]&
{\Oplus_{s+t=n+1}B^s\otimes\tilde B^t}
\ar[r] & 0\\
\\
0\ar[r]&{\Oplus_{s+t=n}Z^s\otimes\tilde B^t\push{B^s\otimes\tilde B^t}
B^s\otimes\tilde Z^t\ }
\ar[0,1]\ar@{>->}[-2,0]&
{E(1)_{-n}(E_{\gamma_n})}
\ar[0,1]\ar[-2,0]
_{E(1)_{-n}({\tiny\hocolim{}}(p^{\zeta_n}_{\gamma_{n}}))}
&
{\Oplus_{s+t=n+1}B^s\otimes\tilde B^t}
\ar@{=}[-2,0]\ar[r] & 0.
} \end{equation}
In particular, $\E (E_{\gamma_n})$ and $\E (E_{\zeta_n})$  are
concentrated in the correct degrees, and
$$E(1)_{-n}(\hocolim{}(p^{\zeta_n}_{\gamma_{n}}))$$ 
is a monomorphism.
This completes the proof of Proposition \ref{inj}.
\end{Pf}
Note that also 
$$\E (E_{\beta_n}) = 
\E (\Oplus_{s+t=n} A_{\beta_s}\wedge\tilde A_{\gamma_t}) = 
\Oplus_{s+t=n} B^s\otimes\tilde B^t$$ is concentrated in the
same degrees.
We need one more step to obtain an object of $\mathcal L$ with the
correct homotopy colimit.
\begin{Rem}
The assumptions of the following proposition are
superfluous. We only state them to simplify the proof.
\end{Rem}
\begin{Prop}\label{i}
Let $E\in\Ho (\EM^{D_N})$ be such that for all $n\in\zn$ and all 
$\alpha\in\{\beta, \gamma, \zeta\}$ the object $\E (E_{\alpha_n})$ is
concentrated in degrees congruent to $-n$ modulo $2p-2$. Let
$i\negmedspace : C_N\to D_N$ send
$\beta_n$ to $\gamma_n$ and $\zeta_n$ to $\zeta_n$.
Then there is an isomorphism, natural in $E$,
$$\hocolim{D_N}E \cong \hocolim{C_N}i^*E.$$
\end{Prop}
\begin{Cor}
Let $E$ be as in Notation \ref{ENot}. Then $i^*E$ is an object of
$\mathcal L$ with
$$\hocolim{C_N}i^*E \cong \Rek (C)\wedge\Rek (\tilde C).$$
Moreover this isomorphism is natural in $C$ and $\tilde C$.
\end{Cor}
\begin{Pf}{ (of Proposition \ref{i})}
We show that the counit of the adjunction 
$$\map{\epsilon_E}{\holkan ii^*E}{E}$$
induces an isomorphism of homotopy colimits.
Since $\E(-)$ is faithful on $\mathcal S_{E(1)}$, it is enough to show
that 
$$\E(\hocolim{{D_N}}\epsilon_E)$$
is an isomorphism. We compute this using spectral sequence (\ref{SF}),
or rather \cite[1.4.35]{Franke:95}. We actually need its construction
and take the notation for it from \cite{Franke:95}. In particular
\begin{equation}\label{resolution}
E\leftarrow\tilde PE\leftarrow\tilde RE\leftarrow\Omega E
\end{equation}
is an exact triangle with $\E$-homology the first step of a
$\colim{}$-acyclic resolution of $\E(E)$. More precisely, it is the resolution
\begin{center}$$
\haekelmuster{0}{0}{Z}
{0}{G}{B}{E(1)_{-n}(E)=}
\hspace{1cm}
\haekelmuster{B}{G\oplus B}{A_1}
{B}{A_2}{B}{E(1)_{-n}(\PP E)=},
$$
\end{center}
where $A_1 = Z\oplus G\oplus B$ and $A_2 = G\oplus B.$
\begin{center}$$
\haekelmuster{B}{G\oplus B}{A_3}
{B}{A_4}{0}
{E(1)_{-n}(\tilde R E) =},$$
\end{center}
where 
$$A_3 = \ker(i_B^Z+i_G^Z+1_Z) \cong B\oplus G$$ and 
$$A_4 = \ker(i_B^G+1_G) \cong B.$$ 
Note that
$$E(1)_{-n}\tilde RE$$ is already $\colim{}$-acyclic, so that the
construction stops here. 
Note also that 
the restrictions to $i(C_N)$ form a $\colim{}$-acyclic
resolution\footnote{The $\colim{1}$'s of diagrams of the shape $C_N$
are computed as equalizers in a way similar to Lemma \ref{hilfsatz3}.}
of $\E(i^*E)$.
The functor $\mathop{(Ho)LKan}_i$ simply adds a bottom row filled with
zeros\footnote{In Franke's setup, the statement about $\holkan i$
follows from the statement about $\lkan i$ using spectral sequence
(\ref{SF}). With model categories, it follows because $\lkan i$
preserves weak equivalences.}
and therefore preserves the property of a diagram to have
$\colim{}$-acyclic $\E$-homology.
Thus, for computing the homotopy colimit of
$$\holkan ii^*E,$$
we may replace (\ref{resolution}) by
$$\holkan ii^*E \leftarrow \holkan ii^*\PP E \leftarrow \holkan ii^*\tilde R E.
$$
Then $\epsilon_E$, $\epsilon_{\tilde PE}$ and $\epsilon_{\tilde RE}$
induce a map of exact couples. Since the image of $i$ is cofinal in
$D_N$, this map becomes an isomorphism in the $E_2$-term.
\end{Pf}
\subsection{Computing $Q(i^*E)$}\label{secondstep}
We have found an object of $\mathcal L$, whose homotopy colimit equals
$\Rek(C)\wedge\Rek(\tilde C)$, namely 
$$i^*E = i^*\holkan{pr}(A\wedge\tilde A).$$
In order to compute
$$
\Rek^{-1}(\Rek(C)\wedge\Rek(\tilde C)),
$$ 
we apply $Q$ to $i^*E$ as explained in Construction \ref{QConst} 
to $i^*E$.
There are two steps: In \ref{objects} we determine the objects, in
\ref{thedifferentials} the differentials of the complex $Q(i^*E)$.
\subsubsection{The objects}\label{objects}
\begin{Not} Write
\begin{eqnarray*}
(\b\to\z) & := & A\at{\beta_{s+1}\leq\zeta_s} \\
\c & := & \cone(\b\to\z), \\
\end{eqnarray*}
and similarly 
$\tilde{\mathcal{B}}^{t+1}$, $\tilde{\mathcal Z}^t$ and
$\tilde{\mathcal C}^t$.
\end{Not}
\begin{Prop}\label{CC}
The cone of a ``diagonal" edge $(\beta_{n+1}\leq\zeta_n)$ of $i^*E$ is
isomorphic to
$$\bigoplus_{s+t=n}\scc,$$
and this isomorphism is natural in $A$ and $\tilde A$.
\end{Prop}
\begin{Pf}{}
In order to apply the cone, we need to compute
$(i^*E)\at{\beta_{n+1}\leq\zeta_n}$ 
as object of $\Ho(\EM^\Ii)$.
By Corollary \ref{Variation}
\begin{equation}\label{schraegekanten1}
(i^*E)\at{\beta_{n+1}\leq\zeta_n}\cong\holkan{p_B}j^*_B(\AAA),
\end{equation}
where $B$, $p_B$ and $j_B$ are as in Notation \ref{BNot}:
{\center 
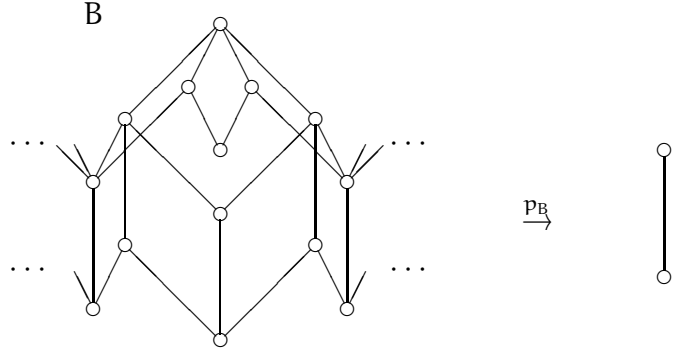
$$
{{\Setunitlength
\begin{picture}(23,11)
\def\Kr{\circle{.4}}
\put(3,11){\makebox(0,0)[cb]{$B$}}
\put(3,2){\Kr}
\put(4,4){\Kr}
\put(7,1){\Kr}
\put(10,4){\Kr}
\put(11,2){\Kr}
\put(7,5){\Kr}
\put(3,6){\Kr}
\put(4,8){\Kr}
\put(10,8){\Kr}
\put(11,6){\Kr}
\put(7,7){\Kr}
\put(7,11){\Kr}
\put(6,9){\Kr}
\put(8,9){\Kr}
\put(21,3){\Kr}
\put(21,7){\Kr}
\put(3,2.2){\line(0,1){3.6}}
\put(4,4.2){\line(0,1){3.6}}
\put(7,1.2){\line(0,1){3.6}}
\put(10,4.2){\line(0,1){3.6}}
\put(11,2.2){\line(0,1){3.6}}
\put(3.0894,2.1789){\line(1,2){.8211}}
\put(3.0894,6.1789){\line(1,2){.8211}}
\put(7.0894,7.1789){\line(1,2){.8211}}
\put(6.0894,9.1789){\line(1,2){.8211}}
\put(10.9106,2.1789){\line(-1,2){.8211}}
\put(10.9106,6.1789){\line(-1,2){.8211}}
\put(6.9106,7.1789){\line(-1,2){.8211}}
\put(7.9106,9.1789){\line(-1,2){.8211}}
\put(11.0894,2.1789){\line(1,2){.5}}
\put(11.0894,6.1789){\line(1,2){.5}}
\put(2.9106,2.1789){\line(-1,2){.5}}
\put(2.9106,6.1789){\line(-1,2){.5}}
\put(7.1414,1.1414){\line(1,1){2.7172}}
\put(6.8586,1.1414){\line(-1,1){2.7172}}
\put(7.1414,5.1414){\line(1,1){2.7172}}
\put(6.8586,5.1414){\line(-1,1){2.7172}}
\put(4.1414,8.1414){\line(1,1){2.7172}}
\put(9.8586,8.1414){\line(-1,1){2.7172}}
\put(3.1414,6.1414){\line(1,1){2.7172}}
\put(10.8586,6.1414){\line(-1,1){2.7172}}
\put(11.1414,6.1414){\line(1,1){1}}
\put(2.8586,6.1414){\line(-1,1){1}}
\put(17,4.5){\makebox(0,0)[cb]{$\stackrel{p_B}{\to}$}}
\put(1,3){\makebox(0,0)[cb]{$\cdots$}}
\put(1,7){\makebox(0,0)[cb]{$\cdots$}}
\put(13,3){\makebox(0,0)[cb]{$\cdots$}}
\put(13,7){\makebox(0,0)[cb]{$\cdots$}}
\put(21,3.2){\line(0,1){3.6}}
\end{picture}
}}
$$}
Let
$$\specialmap{i_B} {B} {\hookrightarrow} {(\CCto\zeta_n)\times\II}$$
be the canonical inclusion, and let
$$\map{pr_\Ii}{(\CCto\zeta_n)\times\II}{\II}$$
be the projection onto the second factor.
Then
\begin{equation}\label{schraegekanten2}
\holkan{p_B}j_B^*\AAA\cong\holkan{pr_\Ii}\holkan{i_B}j_B^*\AAA.
\end{equation}
We need to compute the cone of the right hand side of (\ref{schraegekanten2}).
This right hand side is the homotopy Kan extension along $pr_\Ii$ 
of a diagram which is by (\ref{ecken}) of the form dispayed in Figure 1. 
\begin{figure}[!h]
\begin{center}
$$
\ganzschwarzbutterfly{\szz}{\szb}{\sbz}{\sbb}{\sBz}{\szB}{\sBb}{\sbB}{\sBB}
$$
%
\noindent
\begin{picture}(300,14)
\put(50,14){\makebox(0,0)[cb]{$\uparrow$}}
\put(150,14){\makebox(0,0)[cb]{$\uparrow$}}
\put(250,14){\makebox(0,0)[cb]{$\uparrow$}}
\end{picture}
$$
\npluseinsbutterfly{\bzbbzb}{\sBb}{\sbB}{\star}{\sBz}{\szB}{\sBb}{\sbB}{\sBB.}
$$
\end{center}
Figure 1
\end{figure}
Here $- \push{-} -$ denotes the derived pushout, the vertices in $B$
are black, and those edges that belong to the factor $\II$,
are abbreviated by vertical arrows.
\begin{Rem}
If we would forget some information and think of Figure 1 as a
morphism in
$$\Ho\left(\EM^{\CCto\zeta_n}\right),$$ 
we would obtain the morphisms
$E\at{\beta_{n+1}\leq\zeta_n}$
as homotopy colimit over
$$\CCto\zeta_n$$ of the map in Figure 1.
\end{Rem}
By (\ref{holcone}) together with (\ref{schraegekanten1}) and
(\ref{schraegekanten2}) 
$$\cone((i^*E)\at{\beta_{n+1}\leq\zeta_n})\cong
\hocolim{\CCto\zeta_n}\cone_{\CCto\zeta_n}(\holkan{i_B}j_B^*\AAA).$$
By \cite[1.4.2]{Franke:95},
we know that for a diagram
$X\in\Ho(\MM^{C\times\Ii})$ 
the vertex
$$\left(\cone_C(X)\right)_c$$
is isomorphic to the cone of the corresponding restriction $X\at{c\times\Ii}$.
Therefore,
$$
\cone_{\CCto\zeta_n}\holkan{i_B}j_B^*\AAA
$$
has the form
\begin{center}
\begin{equation}\label{kegelvonschraeg}
\butterfly{\cone\left(\bzbbzb\to\szz\right)}
{\scb}{\sbc}{\Sigma\sbb}{\star}{\star}{\star}{\star}{\star}
\end{equation}
\end{center}
Consider the subposet of $\CCto\zeta_n$
$$\Ww_n := \{(\zeta_s,\zeta_t),(\beta_{s+1},\beta_t)\mid s+t=n\},$$
and let $\jw$ denote its inclusion.
Since $\jw$ has a left adjoint, there is a functor isomorphism
$$\hocolim{\CCto\zeta_n}\cong
\hocolim{\ww}\jw^*.$$
But $\hocolim{\ww}\jw^*$ applied to (\ref{kegelvonschraeg}) becomes
$$
\bigoplus_{s+t=n}\cone\left(\left(\holkan{i_B}j_B^*\AAA\right)
                            _{(\zeta_s,\zeta_t)\times\Ii}\right).
$$
We have already hinted in the pictures, that
$$
\left(\holkan{i_B} j_B^*\AAA\right)_{(\zeta_s,\zeta_t)\times\Ii}
$$
is of the form
$$
\bzbbzb\longrightarrow\szz.
$$
We need to make this more precise.
By Corollary \ref{Kanten}, we have
$$
\left(\holkan{i_B}j_B^*\AAA\right)_{(\zeta_s,\zeta_t)\times\Ii}
\cong\holkan{p_{B'}}j^*_{B'}\AAA,
$$
where
\begin{eqnarray*}
B' & := & B\to (\zeta_s,\zeta_t) \\
p_{B'} & := & p_B\at{B'} \\
j_{B'} & := & j_B\at{B'}. \\
\end{eqnarray*}
This expression can be simplified as follows.
Consider the map of posets
\begin{eqnarray*}
\iC\negmedspace : \II\times\II & \to & B' \\
(0,0) & \mapsto & (\beta_{s+1},\beta_{t+1})\times\{ 0\} \\
(1,0) & \mapsto & (\beta_{s+1},\zeta_t)\times\{ 0\} \\
(0,1) & \mapsto & (\zeta_s,\beta_{t+1})\times\{ 0\} \\
(1,1) & \mapsto & (\zeta_s,\zeta_t).
\end{eqnarray*}
It has a left adjoint $\lC$, satisfying
$$
p_{B'}\circ\iC\circ\lC = p_{B'}.
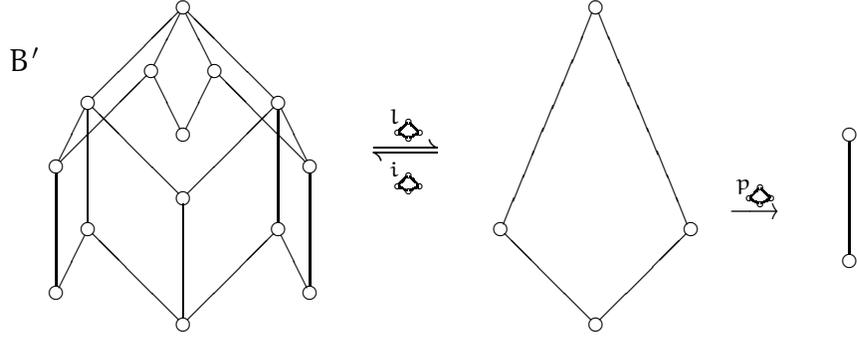
$$
\begin{center}$$
{\Setunitlength
\begin{picture}(29,12)
\def\Kr{\circle{.4}}
\put(3,2){\Kr}
\put(4,4){\Kr}
\put(7,1){\Kr}
\put(10,4){\Kr}
\put(11,2){\Kr}
\put(7,5){\Kr}
\put(3,6){\Kr}
\put(4,8){\Kr}
\put(10,8){\Kr}
\put(11,6){\Kr}
\put(7,7){\Kr}
\put(7,11){\Kr}
\put(6,9){\Kr}
\put(8,9){\Kr}
\put(3,2.2){\line(0,1){3.6}}
\put(4,4.2){\line(0,1){3.6}}
\put(7,1.2){\line(0,1){3.6}}
\put(10,4.2){\line(0,1){3.6}}
\put(11,2.2){\line(0,1){3.6}}
\put(3.0894,2.1789){\line(1,2){.8211}}
\put(3.0894,6.1789){\line(1,2){.8211}}
\put(7.0894,7.1789){\line(1,2){.8211}}
\put(6.0894,9.1789){\line(1,2){.8211}}
\put(10.9106,2.1789){\line(-1,2){.8211}}
\put(10.9106,6.1789){\line(-1,2){.8211}}
\put(6.9106,7.1789){\line(-1,2){.8211}}
\put(7.9106,9.1789){\line(-1,2){.8211}}
\put(7.1414,1.1414){\line(1,1){2.7172}}
\put(6.8586,1.1414){\line(-1,1){2.7172}}
\put(7.1414,5.1414){\line(1,1){2.7172}}
\put(6.8586,5.1414){\line(-1,1){2.7172}}
\put(4.1414,8.1414){\line(1,1){2.7172}}
\put(9.8586,8.1414){\line(-1,1){2.7172}}
\put(3.1414,6.1414){\line(1,1){2.7172}}
\put(10.8586,6.1414){\line(-1,1){2.7172}}
\put(2,9){\makebox(0,0)[cb]{$B'$}}
\put(14,6.5){\makebox(0,0)[ca]
{$\xymatrix{{}\ar @<.2ex> @{{}-^{>}} [r]^{\lC} & {}
\ar @<.2ex> @{{}-^{>}} [l]^{\iC} }$}}
\put(17,4){\Kr}
\put(20,1){\Kr}
\put(23,4){\Kr}
\put(20,11){\Kr}
\put(19.8586,1.1414){\line(-1,1){2.7172}}
\put(20.1414,1.1414){\line(1,1){2.7172}}
\drawline(17.1029,4.1715)(19.8971,10.8285)
\drawline(22.8971,4.1715)(20.1029,10.8285)
\put(25,5){\makebox(0,0)[cc]{$\stackrel{\pC}{\longrightarrow}$}}
\put(28,3){\Kr}
\put(28,7){\Kr}
\put(28,3.2){\line(0,1){3.6}}
\end{picture}}
$$
\end{center}
Therefore, with
$$
\jC := j_{B'}\circ\iC \quad\text{and  }
\pC = P_{B'}\circ \iC 
$$
as in (\ref{box}),
$$
\holkan{p_{B'}}j^*_{B'} \cong \holkan\pC\jC^*\AAA.
$$
By Corollary \ref{smzur},
$$
\jC^*\AAA \cong (A\at{\beta_{s+1}\leq\zeta_s})\wedge
(\tilde A\at{\beta_{t+1}\leq\zeta_t}).
$$
Therefore, by (\ref{lbox}),
$$
\holkan{\pC}\jC^*\AAA\cong (A\at{\beta_{s+1}\leq\zeta_s})\Box
(\tilde A\at{\beta_{t+1}\leq\zeta_t}).
$$
Finally, Proposition \ref{smcone} implies the claim:
$$
\cone\left(\holkan\pC\jC^*\AAA\right)\cong
\cone\left(A\at{\beta_{s+1}\leq\zeta_s}\right)\wedge
\cone\left(\tilde A\at{\beta_{t+1}\leq\zeta_t}\right).
$$
\end{Pf}
\subsubsection{The differentials}\label{thedifferentials}
We have determined the objects of $Q(i^*E)$, and are
ready to compute the differentials.
It turns out that it is enough to consider a  
simple special case:
\begin{Rem}[Franke]\label{FrankeRem}
Let $C^\bullet$ be a cochain complex, and let $s\in\mathbf Z$. 
Consider the map of cochain complexes
\begin{equation}\label{C-C}
\xymatrix{
\cdots \ar[r] & 0 \ar[r]\ar[d] & {C^s}\ar@{=}[r]\ar@{=}[d] & 
{C^s} \ar[r]\ar^{d^s}[d] & 0 \ar[r]\ar[d] & \cdots \\
\cdots \ar[r] & {C^{s-1}} \ar[r] & {C^s} \ar^{d^s}[r] & {C^{s+1}} \ar[r]&
{C^{s+2}} \ar[r] & \cdots,
}
\end{equation}
which we will in the following denote $f_{C,s}$. We write 
$f_{\tilde C,t}$ for the same map for $\tilde C$ and $t$. 
We apply Propositions \ref{inj}, \ref{i} and \ref{CC} to the maps
$f_{C,s}$ and 
$f_{\tilde C,t}$ and obtain
$$\xymatrix{
\cdots\ar[r] & {C^s\otimes\tilde C^t}\ar@{>->}[d]\ar^{?\phantom{hhhhhh}}[r] &
{C^s\otimes\tilde C^t\oplus C^s\otimes\tilde C^t}\ar[r]
\ar^{(d^s\otimes 1,1\otimes\tilde d^t)}[d] & \cdots \\
\cdots \ar[r] & \Oplus_{s'+t'=n} C^{s'}\otimes\tilde
C^{t'}\ar^{?\phantom{lh}}[r]  &
\Oplus_{s'+t'=n+1} C^{s'}\otimes\tilde C^{t'}\ar[r] & \cdots.
}$$
The left vertical arrow is the inclusion of the $(s,t)^{th}$
summand, and the horizontal arrows are the differentials we are
looking for. Therefore,
in order to find the differential in the general case, i.e. the lower
horizontal arrow, it is sufficient to consider the much simpler case
that $C$ and $\tilde C$ are both as in the top row of (\ref{C-C}).
Note that $Q^{-1}$ maps a cochain complex of this form into an
object of $\mathcal L$ that looks like
\begin{equation}\label{CIC}
\xymatrix@=2ex
{&\star \ar@{-}[2,0] \ar@{-}[2,2] && {\mathcal C} \ar@{=}[2,0]
\ar@{-}[2,2] && \star \ar@{-}[2,0]\ar@{-}[1,1] &\\     
{\cdots}\ar@{-}[1,1]& && && &\cdots\\ 
&\star && {\mathcal C} && \star, & \\
&\save*[] {\parbox{10mm}{\raggedright\tiny $s^{st}$ spot}}\restore
\ar@{..>}[u]&& && &\\
}\end{equation}
where
$$C^s=E(1)_*(\Sigma{\mathcal C}).$$
\end{Rem}
\begin{Prop}
Let $A$ and $\tilde A\in\Ob\mathcal L$ be two objects of the form (\ref{CIC})
for $s$, ${\mathcal C}$ and $t$, $\tilde {\mathcal C}$. 
Then the map
(\ref{diff}) that induces the $(s+t)^{th}$ differential in $Q(i^*E)$
is
$$
\map{\diag}{\Sigma^2{\mathcal C}\wedge\tilde{\mathcal C}}
{\Sigma^2{\mathcal C}\wedge\tilde{\mathcal C}\oplus
\Sigma^2{\mathcal C}\wedge\tilde{\mathcal C}.}
$$ 
\end{Prop}
\begin{Pf}{}
Let $n=s+t$.
By Corollary \ref{Kanten}, $i^*E$ looks like
\begin{equation}
\xymatrix@=2ex
{&\star \ar@{-}[2,0] \ar@{-}[2,2] && 
{\Sigma{\mathcal C}\wedge\tilde{\mathcal C}} \ar@{=}[2,0] 
&&{{\mathcal C}\wedge\tilde{\mathcal C}} \ar@{=}[2,0]\ar@{-}[2,2]
&& \star \ar@{-}[2,0]\ar@{-}[1,1] &\\     
{\cdots}\ar@{-}[1,1]& && && && &\cdots\\ 
&\star && {\Sigma{\mathcal C}\wedge\tilde{\mathcal C}} 
&&{{\mathcal C}\wedge\tilde{\mathcal C}}
\ar[-2,-2]_0
&&\star & .\\
&\save*[] {\parbox{10mm}{\raggedright\tiny $n^{th}$ spot}}\restore
\ar@{..>}[u]&& && &\\
}\end{equation}

The first morphism of (\ref{diff}) is
$$\Cone(\Cone(i^*E_{\beta_{n+1}\leq\zeta_n})).$$
In our case this is the identity of
$\Sigma^2{\mathcal C}\wedge\tilde{\mathcal C}$.

The second morphism of (\ref{diff}) is the suspension of
$$(i^*E)_{\beta_{n+1}\leq\zeta_{n+1}},$$
and therefore also equal to the identity of
$\Sigma^2{\mathcal C}\wedge\tilde{\mathcal C}$.

The third morphism of (\ref{diff}) is the suspension of the cone inclusion
belonging to
$$(i^*E)_{\beta_{n+2}\leq\zeta_{n+1}}.$$
Corollary \ref{Kanten} implies that
$(i^*E)_{\beta_{n+2}\leq\zeta_{n+1}}$ is the
``equatorial embedding" of
$\Sigma{\mathcal C}\wedge\tilde{\mathcal C}$
into its suspension, i.e. it is the left homotopy Kan
extension along
$$\VVv\times\II\to\II$$
of
\begin{equation}\label{Aequator}
\xymatrix@=.8ex{ 
{\mathcal C}\wedge\tilde{\mathcal C} &
& {\mathcal C}\wedge\tilde{\mathcal C} && 
\star && \star\\
&&\ar[0,2]&&\\
&{\mathcal C}\wedge\tilde{\mathcal C}\ar@{=}[-2,-1]\ar@{=}[-2,1]&&&&
{\mathcal C}\wedge\tilde{\mathcal C}.\ar@{-}[-2,-1]\ar@{-}[-2,1]\\}
\end{equation}
The claim now follows from the below lemma about the cone inclusion of
the equatorial embedding.
\end{Pf}
\begin{Lem}
The cone inclusion of the equatorial embedding of $X$ into $\Sigma X$,
defined as in (\ref{Aequator}),
is equal to
$$\map{\diag}{\Sigma X}{\Sigma X\oplus\Sigma X}.$$
\end{Lem}
%

\begin{Pf}{}
Taking cones commutes with homotopy colimits. Therefore the morphism
that we want to identify is equal to the homotopy colimit over
$\VVv$ of the lower right hand arrow of the diagram
\begin{equation}\label{4V}
\xymatrix@=.8ex{
&&&&\star&&\star&&&&\\
&&&&&&&\ar[1,2]\\
&\ar[-1,2]&&&&\star\ar@{-}[-2,-1]\ar@{-}[-2,1]&&&&&\\
X && X && \star && \star&&{\Sigma X}&&{\Sigma X}\\
&&\ar[0,2]&&&&\ar[0,2]&&&&\\
&X\ar@{=}[-2,-1]\ar@{=}[-2,1]&&&&
X\ar@{-}[-2,-1]\ar@{-}[-2,1]&&&&\star.\ar@{-}[-2,-1]\ar@{-}[-2,1]\\}
\end{equation}
This arrow is a baby-phantom: restricted to any vertex of $\VVv$, it
becomes zero, yet its homotopy colimit is not equal to zero.
For simplicial sets, and therefore for finite spectra, the claim of
the lemma is well known. It can for example be seen by
looking at explicit cofibrant replacements in the last diagram.
But this special case implies the general case:
In \cite[Cor.1.6.1.]{Franke:95}, Franke defines a family of bi-functors
$$\map{-\wedge-}{\mathcal K_C\times\mathcal S_C^{fin}}{\mathcal
K_{C\times D}}$$ 
that commutes with homotopy colimits and cones, and is associative.
It has the properties that
$-\wedge\mathbb S^0$ 
is the identity,
$-\wedge\mathbb S^1$
is the suspension, and that
$-\wedge\star$
is zero.
Here $\mathcal K_D$ is an arbitrary system of triangulated diagram
categories, in our case
$$\mathcal K_D = \Ho (\MM^D_{S_{E(1)}}),$$
and $S_C^{fin}$ denotes the homotopy category of $C$-shaped diagrams
of finite spectra. In order to obtain the two diagrams
(\ref{Aequator}) and (\ref{4V}) for arbitrary $X$, we take the
corresponding diagrams for $\mathbb S^0$ and smash it with $X$.
Since the claim of the lemma is true for $\mathbb S^0$, and by the
properties of Franke's smash product, it follows that it is also true
for $X$.
\end{Pf}
Together with Remark \ref{FrankeRem} this completes the proof of
Theorem \ref{Maintheorem}.
$$
\Sigma^2{\mathcal C}\wedge\tilde{\mathcal C}\cong
\Sigma{\mathcal C}\wedge\Sigma\tilde{\mathcal C}.
$$
%
%
%
%
\bibliographystyle{alpha}
\bibliography{diplom.english.bib}
\end {document}